\documentstyle[12pt]{article}
\input{amssym.def}
\input epsf

\newtheorem{thm}{Theorem}[section]
\newtheorem{lemma}[thm]{Lemma}
\newtheorem{prop}[thm]{Proposition}
\newtheorem{conj}[thm]{Conjecture}
\newtheorem{defn}[thm]{Definition}

\newtheorem{cor}[thm]{Corollary}
\newtheorem{rmk}[thm]{Remark}

\newcommand{\Cone}{{\rm{Cone}}}

\def\index{{\rm{index}}}

\def\ker{\mathop{\rm Ker}\nolimits}

\def\supp{{\mbox{\rm supp}}}
\def\smax{{\mbox{\rm smax}}}

\newcommand{\C}{{\Bbb C}}

\newcommand{\R}{{\Bbb R}}

\newcommand{\Rgeq}{\R_{\geq 0}}

\newcommand{\PP}{{\Bbb P}}

\newcommand{\HH}{{\Bbb H}}
\newcommand{\cS}{{\cal S}}
\newcommand{\cZ}{{\cal Z}}

\newcommand{\cI}{{\cal I}}
\newcommand{\cJ}{{\cal J}}
\newcommand{\cK}{{\cal K}}
\newcommand{\cL}{{\cal L}}

\newcommand{\cP}{{\cal P}}

\newcommand{\cU}{{\cal U}}

\newcommand{\cD}{{\cal D}}

\begin{document}
\bibliographystyle{plain}

\title{Dimensions of the Ascending and Descending Sets \\
       in Complex Stratified Morse Theory}

\author{Mikhail Grinberg}

\maketitle

\begin{abstract}
This paper is a sequel to \cite{Gr}.  We present a new construction of
gradient-like vector fields in the setting of Morse theory on a complex
analytic stratification.  We prove that the ascending and descending
sets for these vector fields possess cell decompositions satisfying
the dimension bounds conjectured by M. Goresky and R. MacPherson 
in \cite{GM2}.  The vector fields constructed in \cite{Gr} satisfied the
same dimension bounds only ``up to fuzz.''  The new construction has
a number of other advantages.  In particular, it is closer to ``metric
gradient'' intuition.  The idea is to relate a gradient-like vector field in
the neighborhood of a point stratum $\{ p \}$ to a gradient-like vector
field for the distance to $p$ on the complex link of $\{ p \}$.  Similar
results by C.-H. Cho and G. Marelli have recently appeared in
\cite{CM}.  
\end{abstract}

\tableofcontents

\section{Introduction}

\subsection{Main Results}

This paper is a sequel to \cite{Gr}.  It concerns gradient-like vector
fields on complex analytic Whitney stratifications.  Some technical
notions of Whitney stratification theory, such as control data, are
defined somewhat differently by different authors.  We will be using
the versions from \cite{Gr}, with one correction (see Remark
\ref{cdatrmk}), throughout this paper, and will summarize the basic
technical definitions (control data, controlled and weakly controlled
vector fields, etc.) in Section 2.2.

Let $(X, \cS)$ be a Whitney stratified $C^\infty$ manifold.  Given a
$C^\infty$ function $f : X \to \R$, we denote by $\Sigma_f \subset X$
the set of stratified critical points of $f$ (see Definition
\ref{basicsmt}).  The following definitions introduce the main
characters of our story (cf. \cite[Definitions 2.12-13]{Gr}).

\begin{defn}\label{gradlike}
Let $f : X \to \R$ be a $C^\infty$ function.  An $\cS$-preserving
$\nabla f$-like vector field $V$ on an open subset $\cU \subset X$ 
is a weakly controlled vector field, compatible with some system
of control data on $(X, \cS)$, satisfying:

(a)  $V_p = 0$ for all $p \in \Sigma_f \cap \cU$;

(b)  $V_x \, f > 0$ for all $x \in \cU \setminus \Sigma_f$.
\end{defn}

Given a weakly controlled vector field $V$ on an open set 
$\cU \subset X$ and a number $t \in \R$, we write
$\psi_{V,t} : D(\psi_{V,t}) \to \cU$ for the time-$t$ flow of $V$.
Here $D(\psi_{V,t}) \subset \cU$ is the largest open set on which
this flow can be defined (see \cite[Proposition 2.11]{Gr}).

\begin{defn}\label{andmfds}
Let $f : X \to \R$ be a $C^\infty$ function, let $\cU \subset X$
be an open subset, let $V$ be a $\nabla f$-like vector field
on $\cU$, and let $p \in \Sigma_f \cap \cU$.  We define the
descending set $M^-_V (p)$ to be the set of all $x \in \cU$
such that $x \in D(\psi_{V,t})$ for all $t \geq 0$ and
$$\lim_{t \to +\infty} \psi_{V,t} (x) = p.$$
The ascending set $M^+_V (p)$ is defined similarly.
\end{defn}

In the situation of Definition \ref{andmfds}, given an
$\epsilon \in \R$, define
$$L^\pm_{V,\epsilon} (p) = \{ x \in M^\pm_V (p) \; | \; 
f(x) - f(p) = \epsilon \}.$$

\begin{defn}\label{andlink}
The set $L^-_{V,\epsilon} (p)$, for $\epsilon < 0$, is called a descending
link of $p$, if for every $\epsilon_1 \in (\epsilon, 0)$, there is a
bijection $L^-_{V,\epsilon} (p) \to L^-_{V,\epsilon_1} (p)$ given by
identifying points lying on the same trajectory of $V$.  An ascending link
$L^+_{V,\epsilon} (p)$, for $\epsilon > 0$, is defined similarly.
\end{defn}

The existence of ascending and descending links is established by the
following proposition.

\begin{prop}\label{andlinkexist}
In the situation of Definition \ref{andmfds}, assume that $p$
is an isolated (for example, Morse) critical point of $f$.
Then we have the following.

(i)   There exists an $\epsilon_0 > 0$, such that $L^-_{V,\epsilon} (p)$
      is a descending link for every $\epsilon \in (-\epsilon_0, 0)$ and
      $L^+_{V,\epsilon} (p)$ is an ascending link for every
      $\epsilon \in (0, \epsilon_0)$.

(ii)  The ascending and descending links $L^\pm_{V,\epsilon} (p)$
      are compact.

(iii) For every descending link $L^-_{V,\epsilon} (p)$ and every
      $\epsilon_1 \in (\epsilon, 0)$, the bijection $L^-_{V,\epsilon} (p)
      \to L^-_{V,\epsilon_1} (p)$ of Definition \ref{andlink} is a
      homeomorphism; and similarly for ascending links.

\end{prop}

\noindent
{\bf Proof:}  This is an exercise in general topology, using the 
continuity of the time-$t$ flow $\psi_{V,t} : D(\psi_{V,t}) \to \cU$
(see \cite[Proposition 2.11]{Gr}).
\hfill$\Box$

\vspace{.125in}

The following definition can be compared to \cite[Definition 6.1]{Gr}.

\begin{defn}\label{wsss}
(i) A weakly stratified subset $A \subset X$ is a closed subset,
presented as a finite disjoint union
$\displaystyle A = \bigcup_{i=0}^{N} A_i$ so that:

(a) each $A_i$ is a locally closed smooth submanifold of one of the
     strata of $\cS$;

(b) the partial union $\displaystyle A_{\leq n} = \bigcup_{i=0}^{n} A_i$
     is closed for every $n \in \{0, \dots, N\}$.

(ii) A weakly stratified subset $A \subset X$ is called a cellular subset
if every $A_i$ is diffeomorphic to an open ball.
\end{defn}

Given a weakly stratified subset $A = \bigcup A_i \subset X$, we will
refer to the pieces $A_i$ as the strata of $A$.  When $A$ is a cellular,
we will refer to the $A_i$ as the cells of $A$.  Note that an ordering
of the strata is part of the structure of a weakly stratified subset.
The definition of a weakly stratified subset may seem too weak to have
any useful consequences.  Nevertheless, it suffices to prove a basic
stratified general position result (Lemma \ref{genpos}) which can be
used in applications to self-indexing (see Section 1.3).  Our main result
is the following.

\begin{thm}\label{geom}
Let $X$ be a nonsingular complex analytic variety.  Let $\cS$ be a
complex analytic Whitney stratification of $X$.  Let $p \in X$ be a point
stratum of $\cS$.  Let $f : X \to \R$ be a smooth function.  Assume that
$p$ is a stratified Morse critical point of $f$.  Then there exist an open
neighborhood $\cU \subset X$ of $p$ and an $\cS$-preserving
$\nabla f$-like vector field $V$ on $\, \cU$, such that the following four
conditions hold.

(i)   The ascending and descending sets $M^\pm_V (p)$ are cellular
       subsets of $\cU$.

(ii)  The flow of $V$ preserves the cells of $M^\pm_V (p)$.

(iii) The ascending and descending links $L^\pm_{V,\epsilon} (p)$ are
      cellular subsets of $\cU$, with cell decompositions given by
      intersecting with the cells of $M^\pm_V (p)$.

(iv) For every stratum $B \in \cS$, we have:
      $$\dim_\R M^\pm_V (p) \cap B \leq \dim_\C B.$$
\end{thm}

For background on Theorem \ref{geom} the reader is referred to
the introduction of \cite{Gr}.  Theorem \ref{geom} addresses a
conjecture of Goresky and MacPherson stated in
\cite[Part II, \S 6.6]{GM2} and paraphrased as
\cite[Conjecture 1.4]{Gr}.  It provides an alternative to our
earlier existence results for gradient-like vector fields
near a point stratum \cite[Theorems 1.5, 6.2]{Gr}.  More
precisely, Theorem \ref{geom} can be seen as ``removing the fuzz"
from the statements of those results.  However, Theorem \ref{geom}
is not logically stronger than \cite[Theorems 1.5, 6.2]{Gr}.
The main weakness of Theorem \ref{geom} is that it provides
ascending and descending sets which are merely cellular subsets
of $X$, with no claim about how the cells attach to each other.  One
clear advantage of Theorem \ref{geom} is its simplicity of statement.
Another is that it applies to complex analytic stratifications, while all
stratifications in \cite{Gr} were assumed to be complex algebraic.

The original Goresky-MacPherson conjecture (GMC) differs from
Theorem \ref{geom} in two ways.  First, GMC asserts the
existence of a controlled vector field $V$ on the punctured
neighborhood $\cU \setminus \{ p \}$.  The reason for excluding the
point $\{ p \}$ is that a controlled vector field on $\cU$ can not
have non-trivial trajectories approaching $p$ as time tends to
infinity.  By contrast, Theorem \ref{geom} provides a weakly
controlled vector field $V$ defined on $\cU$.  This seems to be
a technical distinction.  The author is not aware of any strong
reason to be interested in controlled vector fields on
$\cU \setminus \{ p \}$ rather than weakly controlled vector fields
on $\cU$.  The second distinction is more important.  GMC
requires that the ascending and descending sets $M^\pm_V (p)$ be
Whitney stratified, while Theorem \ref{geom} provides
$M^\pm_V (p)$ which are only cellular.  By \cite{Go}, every Whitney
stratification can be refined to a triangulation.  Therefore, every
Whitney stratified subset can be presented as a cellular subset.
But the converse is certainly false.  In this regard, GMC is much
stronger than Theorem \ref{geom}.

Theorem \ref{geom} has a natural generalization to the case of
a Morse critical point which is not a point stratum.  To state it,
we must recall the definition of Morse index in complex stratified 
Morse theory (cf. \cite[Definition 1.3]{Gr}).

\begin{defn}\label{csmtindex}
Let $(X,\cS)$ and $f: X \to \R$ be as in Theorem \ref{geom}.  Let
$p \in \Sigma_f$ be a Morse critical point lying in a stratum
$A \in \cS$.   We define
$$\index_f (p) = \index_{f|_A} (p) - \dim_\C A.$$
\end{defn}

\begin{thm}\label{geom2}
Let $(X,\cS)$, $f: X \to \R$, $p \in \Sigma_f$, and $A \in \cS$
be as in Definition \ref{csmtindex}.  Then there exist an open
neighborhood $\cU \subset X$ of $p$ and an $\cS$-preserving
$\nabla f$-like vector field $V$ on $\cU$, such that conditions
(i)-(iii) as in Theorem \ref{geom} and condition (iv) below hold.

(iv) For every stratum $B \in \cS$, we have:
$$\dim_\R M^-_V (p) \cap B \leq \dim_\C B + \index_f (p),$$
$$\dim_\R M^+_V (p) \cap B \leq \dim_\C B - \index_f (p).$$
\end{thm}

The present author's interest in the Goresky-MacPherson conjecture
was re-awakened by an earlier version of the preprint \cite{CM},
which claims results very similar to Theorems \ref{geom}, \ref{geom2}
(see \cite[Theorem 1.1, Corollary 1.2]{CM}).  A brief correspondence
with Cho and Marelli led to some revisions of their preprint.  In the
meantime, the author discovered the construction of gradient-like
vector fields presented in this paper.  Not having studied the proofs in
\cite{CM}, it appears that the construction of this paper is substantially
different from the construction of Cho and Marelli.  Furthermore, it
seems likely that applications of Theorems \ref{geom}, \ref{geom2},
or of the results of \cite{CM}, will utilize the specific constructions of
gradient-like vector fields used in the proofs, rather than the existence
statements alone (see Section 1.3 for some speculations about
possible applications).  For these reasons, the author hopes that this
paper will not be entirely superfluous.

\subsection{Conjectured Refinements}

Theorem \ref{geom2} can be viewed as an analogue of the classical fact that
the gradient flow of a Morse function $f$ on a compact Riemannian manifold
$X$ gives rise to a cell decomposition of $X$, whose cells are parameterized
by the critical points of $f$ (see \cite{Th1}).  This geometric result admits
a homotopy level refinement asserting that $X$ is homotopy equivalent to a
CW-complex, whose cells are likewise parameterized by the critical points of
$f$ (see \cite[Theorem 3.5]{Mi}).  In fact, it is this refinement that makes
the Morse-theoretic cell decomposition so useful in smooth manifold topology.
In this section, we state as conjectures two homotopy level refinements of
Theorem \ref{geom2} (Conjectures \ref{homotopy}, \ref{homotopy2}).  These
conjectures can be viewed as analogues of \cite[Theorem 3.5]{Mi}.  It seems
likely that they can be established by adapting the classical proof of
\cite[Theorem 3.5]{Mi} to  the inductive combinatorics of cells in
$M^\pm_V (p)$ described by Proposition \ref{celldecomp2}.

We begin with a point of notation.  In the situation of Proposition
\ref{andlinkexist}, the descending links $L^-_{V,\epsilon} (p)$, for
different $\epsilon < 0$, are naturally in bijection with each other.
These bijections are provided by Definition \ref{andlink}.  Moreover, in
the situation of Theorem \ref{geom2}, these bijections preserve all the
structures we are interested in: they are cell-preserving homeomorphisms
which restrict to diffeomorphisms of the individual cells.  For this reason,
we will sometimes drop the $\epsilon$ from the notation, writing $L^-_V (p) =
L^-_{V,\epsilon} (p)$ for the descending link.  Similarly, we will write
$L^+_V (p) = L^+_{V,\epsilon} (p)$ for the ascending link.

Theorem \ref{geom2} provides a local topological picture of $M^-_V (p)$
as a cone over the compact cellular subset $L^-_V (p) \subset \cU$
(and similarly for $M^+_V (p)$).  From a homotopy theory point of view,
it is tempting to ask whether the cell decomposition of $L^-_V (p)$
makes it into a CW-complex.  This is more than we can assert.  However,
we conjecture that there exists a CW-complex $K^-_V (p)$, whose cells
are in one-to-one correspondence with the cells of $L^-_V (p)$, and a
homotopy equivalence $h^-_V (p) : K^-_V (p) \to L^-_V (p)$ respecting
the stratification $\cS$.

To make the above precise, we need the language of $\cJ$-filtered spaces
(see \cite[Part III, \S 2.1]{GM2}).  Let $\cJ$ be a partially ordered
set with a unique maximal element $J_0 \in \cJ$.

\begin{defn}\label{jfiltered}
A $\cJ$-filtration of a topological space $X$ is a collection of closed subsets
$\{ X_J \}_{J \in \cJ}$ such that $X_{J_0} = X$ and $J_1 < J_2 \Rightarrow
X_{J_1} \subset X_{J_2}$.  A $\cJ$-filtered map $f : X \to Y$ between two
$\cJ$-filtered spaces is a continuous map such that $f(X_J) \subset Y_J$ for
every $J \in \cJ$.  $\cJ$-filtered homotopies between $\cJ$-filtered maps and
homotopy equivalences between $\cJ$-filtered spaces are defined in the
obvious way.  A pair of $\cJ$-filtered spaces is a pair $(X,Y)$ of topological
spaces plus a $\cJ$-filtration of $X$.   $\cJ$-filtered maps, homotopies, and
homotopy equivalences for pairs are defined in the obvious way.
\end{defn}

In the situation of Theorem \ref{geom2}, let $\cI$ be the set of all closed
unions of strata of $\cS$, partially ordered by inclusion.  Note that every
closed subset of $X$ is naturally an $\cI$-filtered space.  Let
$\cL^\pm_V (p)$ be the set of cells of $L^\pm_V (p)$.

\begin{conj}\label{homotopy}
The statement of Theorem \ref{geom2} can be strengthened to assert the
following.  There exists a CW-complex $K^-_V (p)$, whose set of cells
we denote by $\cK^-_V (p)$, and a bijection $k^- : \cK^-_V (p) \to \cL^-_V (p)$,
such that the following conditions hold.

(i)   For every $C \in \cK^-_V (p)$, we have $\dim_\R C = \dim_\R k^- (C)$.

(ii)  The bijection $k^-$ makes $K^-_V (p)$ into an $\cI$-filtered space.

(iii) There exists a homotopy equivalence of $\cI$-filtered spaces
      $$h^-_V (p) : K^-_V (p) \to L^-_V (p).$$

\noindent
Similarly, there exist a CW-complex $K^+_V (p)$, a dimension-preserving
bijection $k^+ : \cK^+_V (p) \to \cL^+_V (p)$, and a homotopy equivalence
of $\cI$-filtered spaces
$$h^+_V (p) : K^+_V (p) \to L^+_V (p).$$
\end{conj}

Our next conjecture ties in the ascending and descending
links $L^\pm_V (p)$ with a central ingredient of stratified Morse theory:
the local Morse data (see \cite[Part I, \S 3.5]{GM2}).  We recall
the definition of the local Morse data in the situation of Definition
\ref{csmtindex}.  Fix a Riemannian metric $\mu$ on $X$ and a pair of
numbers $0 \ll \epsilon \ll \delta \ll 1$.  Let $B^\mu_\delta (p)
\subset X$ be the closed $\delta$-ball around $p$.  Define
$$D_f(p) = \{ x \in B^\mu_\delta (p) \; | \; f(x) - f(p) \in [-\epsilon, \epsilon] \},$$
$$E^\pm_f(p) = \{ x \in B^\mu_\delta (p) \; | \; f(x) - f(p) = \pm \epsilon \}.$$
The local Morse data for $f$ at $p$ is the pair $(D_f(p),E^-_f(p))$ of
closed subspaces of $X$.  Similarly, the pair $(D_f(p),E^+_f(p))$ is
the local Morse data for $-f$ at $p$.  We suppress from the notation
the inessential dependence of the pairs $(D_f(p),E^\pm_f(p))$ on $\mu$,
$\delta$, $\epsilon$ (see \cite[Part I, Theorem 7.4.1]{GM2}).  In the situation
of  Theorem \ref{geom2},  we have natural inclusion maps:
$$l^\pm_V (p) : L^\pm_V (p) \to E^\pm_f(p).$$

\begin{conj}\label{homotopy2}
The statement of Conjecture \ref{homotopy} can be strengthened to
assert the following.  The maps $l^\pm_V (p) : L^\pm_V (p) \to E^\pm_f(p)$
are homotopy equivalences of $\cI$-filtered spaces.
\end{conj}

Using a natural generalization of \cite[Part I, Theorem 3.12]{GM2}
to describe the $\cI$-filtered homotopy type of the local Morse data,
we obtain the following corollary of Conjecture \ref{homotopy2}.

\begin{cor}\label{homotopy2pairs}
In the situation of Conjecture \ref{homotopy2}, the natural extensions
of $l^\pm_V (p)$ to inclusion maps of pairs:
$$\tilde l^\pm_V (p) : (M^\pm_V (p) \cap D_f (p), L^\pm_V (p)) \to
(D_f (p), E^\pm_f (p)),$$
are homotopy equivalences of pairs of $\cI$-filtered spaces.
\end{cor}

Note that, by combining the dimension bounds of Theorem \ref{geom2}
with the homotopy equivalences of $\cI$-filtered spaces:
$$l^\pm_V (p) \circ h^\pm_V (p) : K^\pm_V (p) \to E^\pm_f(p),$$
provided by Conjectures \ref{homotopy}, \ref{homotopy2}, we obtain
dimension bounds for the $\cI$-filtered homotopy type of $E^\pm_f(p)$.
These dimension bounds are not new.  They are essentially contained,
minus the language of $\cI$-filtered homotopy type, in
\cite[Introduction, \S 1.5]{GM2} and \cite[Part II, \S $1.1^*$]{GM2}.
The latter section also contains references to the original work
of Kato, Karchyauskas, and Hamm proving that a Stein space of
complex dimension $n$ is homotopy equivalent to a CW-complex of
(real) dimension $n$, which is a closely related phenomenon.

\subsection{Some Motivation}

In this section, we briefly indicate four reasons why one might
be interested in the Goresky-MacPherson conjecture (GMC),
Theorems \ref{geom}, \ref{geom2} and Conjectures
\ref{homotopy}, \ref{homotopy2} stated above, and the results
of \cite{Gr}.

The first reason is that vector fields of the kind constructed in this
paper provide a geometric ingredient making some of the main
results of complex stratified Morse theory more concrete and
geometrically apparent.  One example of this is the homotopy
dimension bounds for the local Morse data discussed at the end
of Section 1.2.  Another example is the vanishing of intersection
homology Morse groups outside of a single degree (see
\cite[Part II, Theorem 6.4]{GM2}).  In fact, GMC was introduced
as a source of geometric intuition behind this result.

The second reason is that the results of this paper may enable
some progress towards another conjecture of Goresky and
MacPherson, stated in \cite[Part II, \S $1.1^*$]{GM2}.  Here we
state a variant of this conjecture which is likely to be accessible
using the vector fields provided by the proof of Theorem \ref{geom2}.

\begin{conj}\label{retract}
Let $X \cong \C^d$ be a complex vector space, and let $\cS$ be
a complex analytic Whitney stratification of $X$.  Let $\mu$ be a
Hermitian metric on $X$, let $p \in X$ be the origin, and let
$f : X \to \R$ be the function $f(x) = \mbox{\em dist}_\mu(p,x)^2$. 
Assume that $f$ is Morse.  Pick a regular value $a >0$ of $f$.
Let $X^a = \{ x \in X \; | \; f(x) \leq a \}$ and let $\Sigma^a_f = 
\{ p \in \Sigma_f \; | \;  f(p) < a \}$.  Then there exists a
$\nabla f$-like vector filed $V$ on $X$ with the following
property.  Define
$$X^a_c = \bigcup_{p \in \Sigma^a_f} M^-_V (p).$$
Then $X^a_c$ is a cellular subset of $X$, satisfying
$\dim_\R X^a_c \cap B \leq \dim_\C B$ for every $B \in \cS$.
Moreover, $X^a$ deformation retracts to $X^a_c$ by a stratum
preserving retraction.
\end{conj}

The third reason is that Theorem \ref{geom2} can be used to
prove the existence of self-indexing Morse functions for complex
analytic stratifications.  More precisely, we have the following
generalization of \cite[Theorem 1.6]{Gr} from the complex
algebraic to the complex analytic setting.

\begin{defn}\label{seflind}
Let $X$ be a compact, nonsingular complex analytic variety, and
let $\cS$ be a complex analytic Whitney stratification of $X$.
A Morse function $f : X \to \R$ is called self-indexing if
$f(p) = \index_f (p)$ for every $x \in \Sigma_f$.
\end{defn}

\begin{thm}\label{siexist}
For every pair $(X, \cS)$ as in Definition \ref{seflind}, there
exists a self-indexing Morse function $f : X \to \R$.
\end{thm}

\noindent
{\bf Proof:}  This is a straightforward adaptation of the proof
of \cite[Theorem 1.6]{Gr}, using Theorem \ref{geom2} in place of
\cite[Theorem 6.3]{Gr} and Lemma \ref{genpos} in place of 
\cite[Lemma 6.5]{Gr}.
\hfill$\Box$

\vspace{.125in}

The following lemma is analogous to \cite[Lemma 6.5]{Gr}.
It uses the notion of a time-dependent controlled vector
field which is discussed in \cite[\S 6.2]{Gr}.

\begin{lemma}\label{genpos}
Let $(X, \cS)$ be a Whitney stratified $C^\infty$ manifold with a
fixed system of control data $\cD$.  Let $A, B \subset X$ be two
weakly stratified subsets.  Assume $A \cap B$ is compact and, for
every $S \in \cS$, we have:
$$\dim (A \cap S) + \dim (B \cap S) < \dim S.$$
Then there exists a time-dependent controlled vector field
with compact support $\{ V_t \}_{t \in (0,1)}$ on $X$
compatible with $\cD$, whose time-1 flow satisfies
$\psi_{V,1} (A) \cap B = \emptyset$.
\end{lemma}

\noindent
{\bf Proof:}  This is similar to the proof of \cite[Lemma 6.5]{Gr},
using the linear ordering of the strata of $A$ and $B$ given by
Definition \ref{wsss} instead of the partial ordering by dimension.
\hfill$\Box$

\vspace{.125in}

Even in the complex algebraic case, using Thorem \ref{geom2}
in place of \cite[Theorem 6.3]{Gr} simplifies the proof of
Theorem \ref{siexist} .  This is because the statement of
\cite[Theorem 6.3]{Gr} is more complicated, involving the
choice of an open set informally referred to as ``the fuzz.''
The reader is referred to \cite[\S 1.3]{Gr} and the lecture notes
\cite{MP}, \cite{GrM} for a discussion of the role of self-indexing
Morse functions in the theory of middle perversity perverse
sheaves.

The fourth reason to be interested in the results of this paper is that
the sets $M^\pm_V (p)$  provided by Theorem \ref{geom2} can serve
as building blocks of a useful cycle theory for (the hypercohomology
with coefficients in) middle perversity perverse sheaves.  For
example, in the situation of Conjecture \ref{retract}, let $\cP(X, \cS)$
be the category of middle perversity perverse sheaves on $(X, \cS)$
with coefficients in $\C$.  Let $P \in \cP(X, \cS)$, and let $u \in
\HH^k (X^a, \partial X^a; P)$ be a relative hypercohomology class
with coefficients in $P$.  For $i \in \{ -d, -d+1, \dots, 0 \}$, define
$$\Sigma^a_f [i] = \{ p \in \Sigma^a_f \; | \; \index_f (p) \leq i \},$$
$$X^a_c [i] = \bigcup_{p \in \Sigma^a_f [i]} M^-_V (p).$$
Conjecture \ref{retract} can be strengthened to assert that each
$X^a_c [i] \subset X$ is a cellular subset satisfying
$$\dim_\R X^a_c [i] \cap B \leq \dim_\C B + i,$$ 
for every $B \in \cS$.  In this situation, we can think of $X^a_c [-k]
\subset X^a_c$ as a geometric cycle representing (or supporting)
the class $u$.  It seems plausible that a geometric cycle theory
along these lines can be useful in understanding the structure of
the category $\cP(X, \cS)$.

\subsection{Contents of this Paper}

The rest of this paper is organized as follows.  Section 2 is, in large
part, a summary and paraphrase (and in one instance a correction) of
the definitions and results of \cite[\S\S 2-3]{Gr}.  Section 3 is devoted to
a model example: the Hermitian metric gradient of a linear function on
an affine cone $X \subset \C^n$ over a smooth projective variety $PX
\subset \C\PP^{n-1}$.  It turns out that there is a simple geometric reason
why the ascending and descending sets have half the dimension of $X$
in this case.  The main idea of this paper is to maneuver the general
case of Theorem \ref{geom} to look like this model example.  Section 4
describes an inductive proof of Theorems \ref{geom} and \ref{geom2}
based on this idea.  We have tried to keep the notation parallel
between Section 3 and 4, to emphasize the analogy.

\section{Technical Preliminaries}

Whitney stratification theory has a reputation as a rather
technical subject.  The present author feels that part of
the reason for this is the fact that there is no canonical,
or intrinsic, notion of an automorphism of a Whitney
stratified $C^\infty$ manifold $(X, \cS)$.  Going back to
the work of Thom \cite{Th2} and Mather \cite{Ma}, in order
to construct useful automorphisms, one fixes additional
structure on $(X, \cS)$, called a system of  control data.
The choice of control data is always somewhat arbitrary,
but subsequent constructions must respect it.  If the
setting includes another geometric ingredient, such as
a Morse function, it may be difficult to simultaneously
keep track of the chosen control data and the Morse
function.

In \cite[\S\S 2-3]{Gr} we introduced a technical innovation
which allowed us to simultaneously keep track of a Morse
function $f$ and a system of control data which is
specifically chosen to be ``adapted'' to $f$.  To achieve
this, we used a notion of control data which is more
flexible than the ones in such references as \cite{G-al}
and \cite{dPW}.  Specifically, we allowed quasi-distance
functions modeled on arbitrary quasi-norms rather than
Euclidian or Riemannian distances.  Given this flexibility,
we used the function $f$ itself in constructing one of
the quasi-distance functions making up the control data.
Unfortunately, in preparing this paper, we discovered a
technical error in \cite[\S 3]{Gr}, which is most readily
corrected by adding an additional compatibility condition
in the definition of control data.  This is explained in Remark
\ref{cdatrmk} below.

The contents of this section are as follows.  In Section 2.1,
we recall some standard definitions pertaining to conormal
varieties and stratified Morse functions.  In Section 2.2, we
define control data, controlled and weakly controlled vector
fields, and some related notions.  For the most part, these
definitions are copied form \cite[\S 2]{Gr}.  However, there
is one significant distinction in the definition of control data,
which gives us a chance to correct an error in \cite[\S 3]{Gr} 
(see Remark \ref{cdatrmk}).  In Section 2.3, we introduce the
notion of $f$-adapted control data.  This is a new way to
formalize the technical innovation of \cite[\S\S 2-3]{Gr}.
It will play a key role in the proofs of our main results in
Section 4.

\subsection{Stratified Morse Functions}

Let $(X, \cS)$ be a Whitney stratified $C^\infty$ manifold, and let
$f : X \to \R$ be a smooth function.  The following definitions are
standard in the subject (cf. \cite[Definition 1.2]{Gr})

\begin{defn}\label{basicsmt}
(i)  Let $p \in X$ be a point contained in a stratum $S$.  We say that
$p$ is critical for $f$ ($p \in \Sigma_f$) if it is critical for the
restriction $f|^{}_S$.

(ii)  For every stratum $S \in \cS$, let $\Lambda_S$ be the conormal
bundle $T^*_S X \subset T^*X$, and let $\Lambda = \Lambda_\cS =
\bigcup_S \Lambda_S$.  The set $\Lambda$ is called the conormal
variety to $\cS$.  By Whitney's condition (a), $\Lambda \subset T^*X$
is a closed subset.  Note that $p \in \Sigma_f$ if and only if $d_p f \in
\Lambda$. 

(iii) Let $S \in \cS$ be a stratum.  A covector $\xi \in \Lambda_S$
is said to be generic if it does not annihilate any of the limits
of tangent spaces to strata $T \in \cS$ with $S \subset
\overline{T}$.  The set of all generic $\xi \in \Lambda_S$ is
denoted by $\Lambda^0_S$.  We also write $\Lambda^0 =
\Lambda^0_\cS = \bigcup_S \Lambda^0_S$.

(iv)  Let $p \in \Sigma_f$, and let $S \in \cS$ be the stratum
containing $p$.  We say that $p$ is Morse for $f$ if it is Morse
for the restriction $f|^{}_S$ and $d_p f \in \Lambda^0$.

(v)   We say that $f$ is a Morse function if every $p \in \Sigma_f$
is Morse for $f$.
\end{defn}

\subsection{Control Data, etc.}

We will use the notation $\R_+ = (0, +\infty) \subset \R$ and
$\Rgeq = [0, +\infty) \subset \R$.

\begin{defn}\label{proj}
Let $X$ be a smooth manifold, let $S \subset X$ be a locally
closed smooth submanifold, and let $U \subset X$ be an open
neighborhood of $S$.  A tubular projection $\Pi : U \to S$ is a 
smooth submersion restricting to the identity on $S$.
\end{defn}

\begin{defn}\label{qdprelim}
Let $M$ be a smooth manifold, and let $E \to M$ be a vector bundle
with zero section $Z$.  A quasi-norm on $E$ is a smooth function
$\rho : E \setminus Z \to \R_+$ such that $\rho(\lambda \, e) =
\lambda \, \rho(e)$ for every $e \in E \setminus Z$ and
$\lambda \in \R_+$.
\end{defn}

Assuming $E \setminus Z \neq \emptyset$, a quasi-norm
$\rho : E \setminus Z \to \R_+$ has a unique continuous extension
$\tilde \rho : E \to \Rgeq$, which is equal to zero on $Z$.
The extension $\tilde \rho$ is not differentiable on $Z$.

\begin{defn}\label{qd}
Let $X$ be a smooth manifold, let $S \subset X$ be a locally
closed smooth submanifold, let $U \subset X$ be an open
neighborhood of $S$, and let $\Pi : U \to S$ be a tubular projection.
A quasi-distance function $\rho : U \setminus S \to \R_+$
compatible with $\Pi$ is a smooth function satisfying the
following condition.  There exist a vector bundle $\pi : E \to S$
with zero section $Z$, an open neighborhood $U' \subset E$ of
$Z$, and a diffeomorphism $\phi : U' \to U$, such that
$\phi |^{}_Z = \pi |^{}_Z$,  $(\Pi \circ \phi) |^{}_{U'} = \pi |^{}_{U'}$,
and $\rho \circ \phi : U' \setminus Z \to \R_+$ is the restriction
of a quasi-norm on $E$.
\end{defn}

Let $(X, \cS)$ be a Whitney stratified $C^\infty$ manifold.

\begin{defn}\label{cdat}
A system of control data on $(X, \cS)$ is a collection
$\{ U_S, \Pi_S, \rho_S\}_{S \in \cS}$, where $U_S \supset S$
is an open neighborhood, $\Pi_S : U_S \to S$ is a tubular
projection, and $\rho_S : U_S \setminus S \to \R_+$ is a
quasi-distance function compatible with $\Pi_S$, subject to
the following conditions for every pair $S, T \in \cS$ with
$S \neq T$.

(1)  We have $U_S \cap T = \emptyset$ unless
$S \subset \partial T$.

(2)  Whenever $S \subset \partial T$, the map $(\Pi_S, \rho_S) :
U_S \cap T \to S \times \R_+$ is a smooth submersion.

(3)  Whenever $S \subset \partial T$, we have $\Pi_S \circ
\Pi_T = \Pi_S$ on $U_S \cap U_T$.

(4)  Whenever $S \subset \partial T$, we have
$\rho_S \circ \Pi_T = \rho_S$ on $(U_S \cap U_T) \setminus S$.
\end{defn}

We will often use single letter notation
$\cD = \{ U_S, \Pi_S, \rho_S \}$ for a system of control data.

\begin{rmk}\label{cdatrmk}
{\em
Definition \ref{cdat} is stronger than the corresponding
\cite[Definition 2.5]{Gr} in two ways.  First, we have added
conditions (1) and (2).  By \cite[Lemma 2.4]{Gr}, these
conditions can always be satisfied by shrinking the
neighborhoods $\{ U_S \}$, so adding them explicitly is
not an important distinction.  Second, we have imposed
the requirement (CR) that $\rho_S$ be compatible with
$\Pi_S$ for every $S \in \cS$.  This distinction is more
important.  In fact, failure to impose CR in
\cite[Definition 2.5]{Gr} led us to make a technical error
in \cite[\S 3]{Gr}.  We take this opportunity to correct this
error.

In several places in \cite[\S 3]{Gr} we implicitly made the
following assumption.  Let $\cD = \{ U_S, \Pi_S, \rho_S\}$
be a system of control data.  Let $A \in \cS$, $a \in A$, and
$N = \Pi_A^{-1} (a)$.  Then $\cD$ restricts to a system of
control data on some open neighborhood $U_a \subset N$
of $a$.  In part, this assumes that, for some $U_a$, the
restriction $\rho^{}_A |^{}_{U_a \setminus \{a\}}$ is a
quasi-distance function.  This is trivially true for control data
in the sense of Definition \ref{cdat}; but it seems quite subtle
and possibly false without CR.  Fortunately, there is no need
to make this assumption.

The most direct fix is to strengthen the definition of control
data by adding the requirement CR.  This necessitates two
other changes.  First, the first sentence of the proof of
\cite[Theorem 3.1]{Gr} must be extended to assume that the
map $p : \cU \to S$ is the restriction of the bundle projection
$X \to M \cong S$.  No other changes are needed in the proof
of \cite[Theorem 3.1]{Gr}, and we will make reference to this
proof in Section 2.3.  Second, we must rephrase the statement
and the proof of \cite[Proposition 3.9]{Gr}, which is an
intermediate step for passing between the local and global
topological stability results \cite[Theorems 3.8, 3.10]{Gr}.
We  sketch the necessary modification briefly. 

We say that a quasi-distance function $\rho : U_S \to S$ is
universally compatible, if it is compatible with every tubular
projection $\Pi : U_S \to S$, after shrinking the neighborhood
$U_S$.  The statement of \cite[Proposition 3.9]{Gr} must be
modified to assert that the control data $\{ U_*, \Pi_*, \rho_*\}$
of \cite[Theorem 3.8]{Gr} can be chosen so that $\rho^{}_{S_W}$
is universally compatible, and to proceed with this assumption.
The proof of \cite[Proposition 3.9]{Gr} then must begin with the
following remark.  The statement of \cite[Theorem 3.1]{Gr} can
be strengthened to assert that the quasi-distance to $B^\circ$,
provided as part of the control data on $U$, is universally
compatible.  To ensure this, it is enough to require, in Step 1
of the proof, that the Euclidean norm of the fiber-wise liner
function $\tilde f : X \to \R$ be the same in every fiber.  It
follows that $\{ U_*, \Pi_*, \rho_*\}$ can be chosen so that
$\rho^{}_{S_W}$ is universally compatible.  The only other
change needed in the proof of \cite[Proposition 3.9]{Gr}
is in the penultimate sentence, where we must ensure that
the restriction $\rho^{}_{S_W} |^{}_{U'}$ is compatible with
the tubular projection $\tilde x \mapsto \mu(\tilde x)$.

The need for this modification of \cite[Proposition 3.9]{Gr} arises
because its proof is the only place in \cite{Gr} where we modify
a tubular projection after the corresponding quasi-distance
function has been fixed.  Unrelated to the above error, the proof
of \cite[Proposition 3.9]{Gr} contains a typo.  Namely, the function
$\theta_{\tilde x}$ must  be defined as distance squared, rather
than simply distance.

The author will evaluate the need for a more formal erratum.
}
\end{rmk}

\begin{defn}\label{compat}
Let $\cD = \{ U_S, \Pi_S, \rho_S \}$ be a system of control data on
$(X, \cS)$, let $\cU \subset X$ be an open subset, let $A$ be a set,
and let $f : \cU \to A$ be a map of sets.  We say that $\cD$ is
$f$-compatible on $\cU$ if, for every $S \in \cS$, there is a
neighborhood $U'_S \subset U_S \cap \cU$ of $S \cap \cU$ such that
$f \circ \Pi_S = f$ on $U'_S$.
\end{defn}

\begin{defn}\label{chom}
Let $(X, \cS)$, $(\hat X, \hat \cS)$ be two Whitney stratified
$C^\infty$ manifolds, with systems of control data $\cD =
\{ U_S, \Pi_S, \rho_S \}$ on $(X, \cS)$ and $\hat \cD =
\{ U_{\hat S}, \Pi_{\hat S}, \rho_{\hat S} \}$ on $(\hat X, 
\hat \cS)$.  A controlled homeomorphism $\phi : X \to \hat X$,
compatible with $\cD$ and $\hat \cD$, is a homeomorphism
which takes strata diffeomorphicly onto strata, establishing a
bijection $S \mapsto \hat S$, and satisfies the following condition. 
For every $S \in \cS$, there is a neighborhood $U'_S \subset U_S$
of $S$, such that $\phi \circ \Pi^{}_S = \Pi_{\hat S} \circ \phi$ on
$U'_S$ and $\rho^{}_S = \rho_{\hat S} \circ \phi$ on $U'_S
\setminus S$.  A controlled homeomorphism between two open
subsets $\, \cU \subset X$ and $\, \hat \cU \subset \hat X$ is
defined similarly (with compatibility conditions imposed in some
neighborhood $U'_S \subset U_S \cap \cU$ of $S \cap \cU$,
for every $S \in \cS$).
\end{defn}

\begin{defn}\label{cvf}
Let $\cD = \{ U_S, \Pi_S, \rho_S\}$ be a system of control data on
$(X, \cS)$ and let $\cU \subset X$ be an open set.  A controlled
vector field $V$ on $\;\cU$ compatible with $\cD$ is a collection
$\{ V_S \}_{S \in \cS}$ of smooth vector fields on the intersections
$S \cap \cU$, satisfying the following condition. For every $S \in
\cS$, there exists a neighborhood $U'_S \subset U_S \cap \cU$
of $S \cap \cU$ such that:

(a) $(\Pi_S)_* V_x = V_{\Pi_S (x)}$ for every $x \in U'_S$;

(b) $V_x \, \rho_S = 0$ for every $x \in U'_S \setminus S$.
\end{defn}

Integrating controlled vector fields is a basic technique for
constructing controlled homeomorphisms, going back to the work
of Thom \cite{Th2} and Mather \cite{Ma}.  The following lemma
(see \cite[Lemma 2.9]{Gr}, \cite[Lemma 4.11]{Sh},
\cite[Theorem 1.1]{dP}) is a basic tool for constructing
controlled vector fields.

\begin{lemma}\label{ecvf}
Let $\cD$ be a system of control data on $(X, \cS)$.  Let $S$ be a
stratum, let $U \subset S$ be open in $S$, and let $V$ be a smooth
vector field on $U$.  Then there exist an open $\cU \subset X$,
with $\cU \cap S = U$, and a controlled vector field $\tilde V$ on
$\cU$ compatible with $\cD$, such that $\tilde V|^{}_U = V$.
Furthermore, $\tilde V$ can be chosen to be continuous as a section
of $T\cU$.
\hfill$\Box$
\end{lemma}

Controlled vector fields are too rigid for discussing ascending
and descending sets.  Indeed, a trajectory of a controlled vector
field can not approach a point on a smaller stratum as time tends
to infinity.  We therefore need the following definition.

\begin{defn}\label{wcvf}
Let $\cD = \{ U_S, \Pi_S, \rho_S\}$ be a system of control data
on $(X, \cS)$ and let $\cU \subset X$ be an open set.  A weakly
controlled vector field $V$ on $\cU$ compatible with $\cD$ is a
collection $\{ V_S \}_{S \in \cS}$ of smooth vector fields on the
intersections $S \cap \cU$, satisfying the following condition.
For every $S \in \cS$, there exist a neighborhood $U'_S \subset
U_S \cap \cU$ of $S \cap \cU$ and a number $k > 0$, such that:

(a) $(\Pi_S)_* V_x = V_{\Pi_S (x)}$ for every $x \in U'_S$;

(b) $|V_x \, \rho_S| < k \cdot \rho_S(x)$ for every
$x \in U'_S \setminus S$.
\end{defn}

By \cite[Proposition 2.11]{Gr}, weakly controlled vector fields
integrate to stratum preserving homeomorphisms.  Unfortunately,
these homeomorphisms do not preserve control data.  For this
reason, weakly controlled vector fields compatible with a given
system of control data do not form a Lie algebra, and we can not
use the flow of one weakly controlled vector to transform 
another.  As a result, we need both controlled and weakly
controlled vector fields.

\subsection{$f$-Adapted Control Data}

We continue with a Whitney stratified $C^\infty$ manifold $(X, \cS)$.
Let $f : X \to \R$ be a smooth function.  Let $A \in \cS$ be a stratum,
and let $K \subset A$ be a compact subset.  In this section, we introduce
the notion of a system of control data on $(X, \cS)$ which is $f$-adapted
near $K$.  This notion will play an important role in the proofs of our
main results in Section 4.  

To begin, we fix a ``smooth maximum'' function
$$\smax : \Rgeq \times \Rgeq \to \Rgeq,$$
satisfying the following properties.

(1)  The restriction $\smax |^{}_{\R_+ \times \R_+}$ is smooth
     and (non-strictly) convex.

(2)  $\smax(x,y) = \max(x,y)$ whenever
     $\max(x,y) \geq 1.1 \cdot \min(x,y)$.

(3)  $\smax(x,y) = \smax(y,x)$ for all $x,y \in \Rgeq$.

(4)  $\smax(a \cdot x, a \cdot y) = a \cdot \smax(x,y)$
     for all $a \in \R_+$ and $x,y \in \Rgeq$. 

\noindent
One can check that (1)-(4) imply:
$$\max(x,y) \leq \smax(x,y) \leq 1.05 \cdot \max(x,y),$$
for all $x,y \in \Rgeq$.

\begin{defn}\label{fadapted}
A system of control data $\cD = \{ U_S, \Pi_S, \rho_S\}$ on $(X, \cS)$
is said to be $f$-adapted near $K$ if there exists an open neighborhood
$\, \cU \subset U_A$ of $K$ such that the following conditions hold.

(1)  Write $\Pi = \Pi_A$, $\rho = \rho_A$, $A_\cU =  A \cap \cU$.  Let
$g : \cU \to \R$ be the function
$$g(x) = f(x) - f(\Pi(x)).$$
Then we have $\Pi (\cU) = A_\cU$, $\Sigma_g = A_\cU$, and
$d_x g \in \Lambda^0_A$ for every $x \in A_\cU$.

(2)  Let $\, \check \cU = \{ x \in \cU \; | \; g(x) = 0 \}$, and let
$\check \Pi = \Pi |_{\check \cU} : \check \cU \to A_\cU$.  Then
$\check \cU$ is a smooth manifold, and $\check \Pi$ is a tubular
projection.  Moreover, there exist a quasi-distance function
$\check \rho : \check \cU \setminus A_\cU \to \R_+$ compatible
with $\check \Pi$, with continuous extension $\hat \rho : \check
\cU \to \Rgeq$, and a tubular projection $\pi : \cU \to \check \cU$,
such that $\check \Pi (\pi (x)) = \Pi(x)$ for every $x \in \cU$, and
$$\rho (x) = \smax (\hat \rho (\pi (x)), \, |g(x)|),$$
for every $x \in \cU \setminus A_\cU$.

(3)  Let
$$\mathring \cU = \{ x \in \cU \; | \; 
\hat \rho (\pi (x)) > 0.8 \cdot |g(x)| \}.$$
Then the map $\Theta = (\Pi, \, g, \, \check \rho \circ \pi) :
\mathring \cU \to A_\cU \times \R \times \R_+$ is a stratified
submersion.

(4)  The system $\cD$ is $\Theta$-compatible on
$\mathring \cU$.
\end{defn}

The usefulness of the notion of $f$-adapted control data comes
primarily from the following two propositions.

\begin{prop}\label{fadaptedexistbasic}
Assume that $f|^{}_A = 0$ and $d_x f \in \Lambda^0_A$ for every
$x \in K$.  Let $\Pi : U \to A$ be a tubular projection.
Then there exists a system of control data $\cD =
\{ U_S, \Pi_S, \rho_S\}$ on $(X, \cS)$ which is $f$-adapted near
$K$ and satisfies $\Pi_A|^{}_\cU = \Pi|^{}_\cU$ for some open
neighborhood $\cU \subset U \cap U_A$ of $K$.
\end{prop}

\noindent
{\bf Proof:}  This is similar to Steps 1-4 of the proof of
\cite[Theorem 3.1]{Gr}.
\hfill$\Box$

\begin{prop}\label{fadaptedflow}
Assume that $f|^{}_A = 0$ and $d_x f \in \Lambda^0_A$ for every
$x \in K$.  Let $\cD$ be a system of control data on $(X, \cS)$ which
is $f$-adapted near $K$.  Let $V$ be a smooth vector field on $A$.
Then there exist an open neighborhood $\, \cU \subset X$ of $K$ and
a controlled vector field $\tilde V$ on $\cU$ compatible with $\cD$,
such that $\tilde V |_{A \cap \, \cU} = V |_{A \cap \, \cU}$ and
$\tilde V_x f = 0$ for every $x \in \cU$.
\end{prop}

\noindent
{\bf Proof:}  This is similar to Step 5 of the proof of
\cite[Theorem 3.1]{Gr}.
\hfill$\Box$

\begin{cor}\label{fadaptedexist}
Let $p \in \Sigma_f \cap A$ be a Morse critical point of $f$.
Then there exists a system of control data $\cD$ on $(X, \cS)$
which is $f$-adapted near $\{ p \}$.
\end{cor}

\noindent
{\bf Proof:}  Pick a tubular projection $\Pi : U \to A$.  Pick
a smooth function $g : X \to \R$ such that $g = f - f \circ \Pi$
in some neighborhood of $p$.  Note that
$d_p g = d_p f \in \Lambda^0_A$.  Use Proposition
\ref{fadaptedexistbasic} to obtain a system of control data
$\cD = \{ U_S, \Pi_S, \rho_S\}$ on $(X, \cS)$ which is
$g$-adapted near $\{ p \}$ and satisfies $\Pi_A|^{}_\cU =
\Pi|^{}_\cU$ for some open neighborhood $\cU \subset
U \cap U_A$ of $p$.  Then $\cD$ is also $f$-adapted near
$\{ p \}$.
\hfill$\Box$

\begin{defn}\label{nsdef}
Let $(X, \cS)$ be a Whitney stratified $C^\infty$ manifold, let
$A \in \cS$ be a stratum, and let $p \in A$.  A normal slice $N$
to $A$ passing through $p$ is a locally closed smooth submanifold
of $X$, such that $p \in N$, $\dim A + \dim N = \dim X$, and $N$
is transverse to the strata of $\cS$.
\end{defn}

A normal slice $N$, as in Definition \ref{nsdef}, is naturally a
Whitney stratified space, with a stratification induced from $\cS$.
The following corollary is a paraphrase of \cite[Corollary 3.2]{Gr}. 
It may be called a ``stratified Morse lemma'' in the sense that it
provides the closest we can come to a local normal form
statement for a stratified Morse function near a critical point.

\begin{cor}\label{fadaptedprod}
In the situation of Corollary \ref{fadaptedexist}, let $\cD =
\{ U_S, \Pi_S, \rho_S\}$ be a system of control data on $(X, \cS)$
which is $f$-adapted near $\{ p \}$.  Let $N = \Pi_A^{-1} (p)$.
Then N is a normal slice to $A$.  Let $\cS_1$ be the stratification
of $N$ induced from $\cS$, let $\cD_1$ be the system of control data
on $(N,\cS_1)$ induced from $\cD$, and let $f_1 = f|^{}_N : N \to \R$.
Then $\cD_1$ is $f_1$-adapted near $\{ p \}$.  Moreover, there exists
an open neighborhood $\cU \subset U_A$ of $p$ with the following
property.  Let $\, \cU_1 = \cU \cap N$ and $\, \cU_2 = \cU \cap A$.
Let $\tilde \cD_1$ be the system of control data on $\cU_1 \times
\cU_2$ induced from $\cD_1$.  Then there exists a controlled
homeomorphism 
$$\phi : \cU_1 \times \cU_2 \to \cU,$$
compatible with $\tilde \cD_1$ and $\cD$, such that for every
$x_1 \in  \cU_1$ and $x_2 \in  \cU_2$, we have:

(i)   $\phi(x_1, p) = x_1$;

(ii)  $\Pi_A \circ \phi \; (x_1, x_2) = x_2$;

(iii) $f \circ \phi \; (x_1, x_2) = f(x_1) + f(x_2) - f(p)$.
\end{cor}

\noindent
{\bf Proof:}  This is similar to the proof of \cite[Corollary 3.2]{Gr}
and to Step 6 in the proof of \cite[Theorem 3.1]{Gr}.  More precisely,
pick a smooth function $g : X \to \R$ such that $g = f - f \circ \Pi_A$
in some neighborhood of $p$.   Let $a = \dim A$, and let
$\{ V_i \}_{i = 1}^a$ be a collection of vector fields on $A$
such that $\{ (V_i)_p \}_{i = 1}^a \subset T_p A$ is a basis.
Proposition \ref{fadaptedflow} provides controlled extensions
$\{ \tilde V_i \}_{i = 1}^a$ of $\{ V_i \}_{i = 1}^a$, defined
in some neighborhood of $p$ and preserving the function $g$.
The homeomorphism $\phi$ is constructed by integrating the
vector fields $\{ \tilde V_i \}_{i = 1}^a$.
\hfill$\Box$

\vspace{.125in}

We have the following remarkable ``softness'' result for
$f$-adapted control data.

\begin{prop}\label{fadaptednsequiv}
Let $(X, \cS)$ be a Whitney stratified $C^\infty$ manifold, and let
$A \in \cS$.  Suppose we have two points $p_0, p_1 \in A$, two
normal slices $N_0 \ni p_0$ and $N_1 \ni p_1$ to $A$, and two
functions $f_0 : N_0 \to \R$ and $f_1 : N_1 \to \R$, with
$f_0(p_0) = f_1(p_1) = 0$.  Assume that the differentials
$d_{p_0} f_0$ and $d_{p_1} f_1$ are both in $\Lambda^0_A$
and, moreover, in the same path-component of $\Lambda^0_A$.
Let $\cS_0$, $\cS_1$ be the stratifications of $N_0$, $N_1$
induced from $\cS$.  Assume that we are given a system of
control data $\cD_0$ on $(N_0, \cS_0)$ which is $f_0$-adapted
near $\{ p_0 \}$, and similarly for $N_1$.  Then there exist open
neighborhoods $\, \cU_0 \subset N_0$ and  $\, \cU_1 \subset
N_1$ of $p_0$ and $p_1$, and a controlled homeomorphism $\phi :
\cU_0 \to \cU_1$ compatible with $\cD_0$ and $\cD_1$, such
that $f_0 |_{\cU_0} = f_1 \circ \phi$.
\end{prop}

\noindent
{\bf Proof:}  The proof consists of two steps.  The first step is to
establish the proposition in the case when $p_0 = p_1$,
$N_0 = N_1$, $f_0 = f_1$, and only the systems of control data
$\cD_0$, $\cD_1$ may differ.  The second step is to deduce the
general case of the proposition.  Both steps are similar to the
proof \cite[Corollary 3.5]{Gr}.  
\hfill$\Box$

\section{A Model Example}

Theorem \ref{geom} is an existence result, and our proof of
it will proceed by construction.  This construction will be
non-canonical, involving a number of arbitrary choices.
However, the idea of the construction comes from a phenomenon
occurring ``in nature.''  Namely the behavior of the metric gradient
of a linear function near an isolated conical singularity.  In this
section, we describe this ``natural phenomenon'' as motivation for
the proof of Theorem \ref{geom} in Section 4.

Let $W \cong \C^n$ be a complex vector space.  Let
$PW \cong \C\PP^{n-1}$ be the associated projective space.
Let $PX \subset PW$ be a smooth algebraic subvariety, and
let $X \subset W$ be the affine cone over $PX$.  Let $p \in W$
be the origin.  Then $X$ has a natural stratification with
two strata $X = X^\circ \cup \{ p \}$, and $W$ has a natural
stratification with three strata $W = W^\circ \cup X^\circ
\cup \{ p \}$.  Fix a Hermitian metric $\mu$ on $W$ and
define $r : W \to \R$ by $r(w) = \mbox{dist}_\mu(p,w)$.  

Pick a linear function $\varphi : W \to \C$ such that the real
part $f = \mbox{Re} (\varphi) : W \to \R$ is a generic covector
at $p$ (i.e., $f \in \Lambda^0_{\{p\}}$ in the notation of
Definition \ref{basicsmt}).  Consider a vector field $V$ on $X$
defined by $V(p) = 0$ and
$$V|_{X^\circ} = r \cdot \nabla_\mu \, f|_{X^\circ}.$$
Then $p$ is a Morse critical point of $f$, and $V$ is the
restriction to $X$ of a  $\nabla f$-like vector field on $W$.  
Define the ascending and descending sets $M^\pm_V (p)
\subset X$ by analogy with Definition \ref{andmfds}.  The main
result of this section (Theorem \ref{model}) is a description
of $M^\pm_V (p)$.
 
Define $W^\pm = \varphi^{-1} (\pm 1) \subset W$.  Consider
the intersections $Y^\pm = W^\pm \cap X$ and the restrictions
$$g^\pm = r|^{}_{Y^\pm}: Y^\pm \to \R.$$
Note that each $Y^\pm$ is a smooth affine subvariety of $W$
and each $g^\pm : Y^\pm \to \R$ is a smooth function.  Consider
the gradient vector fields $G^\pm = \nabla_\mu \, g^\pm$ on
$Y^\pm$.  The flow of $G^\pm$ is integrable for all time.
Define $Y^\pm_c \subset Y^\pm$ to be the union of all bounded
trajectories of $G^\pm$.  More precisely, we let $Y^\pm_c =
\bigcup_\gamma \gamma(\R)$, where $\gamma$ runs over all
parameterized trajectories $\gamma : \R \to Y^\pm$ of $G^\pm$
such that $g^\pm \circ \gamma : \R \to \R$ is bounded.  It is not
hard to check that $Y^\pm_c \subset Y^\pm$ is compact.  Finally,
let $\Cone (Y^\pm_c) \subset X$ be the real cone over $Y^\pm_c$.
More precisely, for $x \in X^\circ$, let $R(p,x) \subset X$ be
the (closed) straight line ray originating from $p$ and passing
through $x$.  Also, let $R^\circ (p,x) = R(p,x) \setminus \{ p \}$.
Then we have:
$$\Cone (Y^\pm_c) = \bigcup_{y \in Y^\pm_c} R(p,y) =
\{ p \} \cup \bigcup_{y \in Y^\pm_c} R^\circ (p,y).$$

\begin{thm}\label{model}
We have $M^\pm_V (p) = \Cone (Y^\pm_c)$.
\end{thm}

Let $\Omega$ be the cone of all Hermitian metrics on $W$, and
let $\Omega^\circ = \Omega^\circ (X,\varphi) \subset \Omega$
be the set of all $\mu \in \Omega$ such that both functions
$g^\pm : Y^\pm \to \R$ are Morse.  

\begin{prop}\label{genmetric}
The set $\Omega^\circ$ is open and dense in $\Omega$.
\end{prop}

\noindent
{\bf Proof:}  Let us focus on the function $g^- : Y^- \to \R$.  Write
$\Sigma_{g^-} \subset Y^-$ for its critical locus.  The conical
property of $X$ and the generic property of $\phi$ imply that, for
every $\mu_0 \in \Omega$, there exist an open neighborhood
$U_{\mu_0} \subset \Omega$ of $\mu_0$ and a compact
$K \subset Y^-$, such that $\Sigma_{g^-} \subset K$ for every
$\mu \in U_{\mu_0}$.  It follows that $\Omega^\circ \subset \Omega$
is open.

Let $\Delta = \ker (\phi) \subset W$.  For a given $\mu \in \Omega$,
let $L = \Delta^\perp \subset W$ be the orthogonal complement to
$\Delta$ relative to $\mu$.  Let $\{ q^- \} = L \cap W^-$.  Then
$g^- : Y^- \to \R$ is Morse if and only if $q^-$ is not a focal
point for the submanifold $Y^- \subset W^-$.  The density of
$\Omega^\circ \subset \Omega$ follows from the nowhere density of
the set $F^- \subset W^-$ of focal points of $Y^-$.  More precisely,
let $\Omega(\Delta)$ be the cone of all Hermitian metrics on $\Delta$,
and let $\delta : \Omega \to \Omega(\Delta)$ be the restriction map.
Then the focal set $F^- \subset W^-$ depends on $\mu \in \Omega$ only
through the image $\delta (\mu)$, and we can move the point $q^-$
away from $F^-$ by perturbing $\mu$ within the fiber $\delta^{-1}
(\delta(\mu))$.
\hfill$\Box$

\vspace{.125in}

Let $d = \dim_\C X$.  We will use the term ``cellular subset of $X$''
to mean ``cellular subset of $W$, contained in $X$.''

\begin{cor}\label{modelcells}
For $\mu \in \Omega^\circ$, the sets $Y^\pm_c$ and
$M^\pm_V (p)$ are cellular subsets of $X$, and we have:
$$\dim_\R Y^\pm_c \leq d - 1 \;\;\; \mbox{and} \;\;\;
\dim_\R M^\pm_V (p) \leq d.$$
\end{cor}

\noindent
{\bf Proof:}  We will only consider $Y^-_c$ and $M^-_V (p)$.
The statement for $Y^+$ and $M^+_V (p)$ is analogous.
Fix a $\mu \in \Omega^\circ$.  The critical locus $\Sigma_{g^-}
\subset Y^-$ is compact (see the first paragraph of the proof
of Proposition \ref{genmetric}).  Therefore, $\Sigma_{g^-}$ is
a finite set.  Order $\Sigma_{g^-}$ by critical value.  More
precisely, fix and ordering $\Sigma_{g^-} = \{ y_0, \dots, y_N \}$,
such that:
\begin{equation}\label{ordercrit}
i \leq j \;\; \Rightarrow \;\;
g^-(y_i) \leq g^-(y_j),
\;\; \mbox{for all} \;\; i,j \in \{ 0, \dots, N \}.
\end{equation}
Note that
\begin{equation}\label{cellstruct}
Y^-_c = \bigcup_{i=0}^N M^-_{G^-} (y_i),
\end{equation}
where each descending set $M^-_{G^-} (y_i)$ is an open cell
of dimension $\index_{g^-} (y_i)$.  For every $n \in \{ 0,
\dots, N \}$, the partial union
$$(Y^-_c)_{\leq n} = \bigcup_{i=0}^n M^-_{G^-} (y_i)$$
is closed because of (\ref{ordercrit}).  This proves that equation
(\ref{cellstruct}) presents $Y^-_c$ as a cellular subset of $X$.
The dimension bound for $Y^-_c$ follows from the inequality:
\begin{equation}\label{indexbound}
\index_{g^-} (y_i) \leq \dim_\C Y^- = d - 1,
\end{equation}
for every $i \in \{ 0, \dots, N \}$ (see \cite[Part I, \S 7]{Mi}).

For every $i \in \{ 0, \dots, N \}$, define
$$M^-_{G^-} (y_i)^\sharp =
\bigcup_{y \in M^-_{G^-} (y_i)} R^\circ (p,y).$$
By Theorem \ref{model}, we have:
\begin{equation}\label{cellstruct2}
M^-_V (p) = \{ p \} \cup
\bigcup_{i=0}^N M^-_{G^-} (y_i)^\sharp.
\end{equation}
Each $M^-_{G^-} (y_i)^\sharp$ is a cell of dimension
$\index_{g^-} (y_i) + 1$.  Equation (\ref{cellstruct2})
presents $M^-_V (p)$ as a cellular subset of $X$.  The
dimension bound for $M^-_V (p)$ follows from
(\ref{indexbound}).
\hfill$\Box$

\vspace{.125in}

We will prove the claim of Theorem \ref{model} about the
descending set $M^-_V (p)$ only.  The claim about
$M^+_V (p)$ is analogous.  Our proof is based on the
following three lemmas.  Let $f^\perp = \mbox{Im}
(\varphi) : W \to \R$ be the imaginary part of $\varphi$.

\begin{lemma}\label{imaginaryintegral}
We have $V_x f^\perp = 0$ for every $x \in X^\circ$.
\end{lemma}

\noindent
{\bf Proof:}
This is a standard consequence of the Hermitian property
of $\mu$.
\hfill$\Box$

\begin{lemma}\label{safeconical}
There exists a $k > 0$, such that
$$r(x) < k \cdot |f(x)|,$$
for every $x \in M^-_V (p)$.
\end{lemma}

\noindent
{\bf Proof:}
Let $\nu : X^\circ \to \R$ be the norm $\nu =
\| \nabla_\mu \, f|_{X^\circ} \|_\mu$.  Since the function
$f$ is linear and the variety $X$ is conical, the function
$\nu$ is constant on each open ray $R^\circ (p,x)$ for
$x \in X^\circ$.  Therefore, $\nu$ attains a minimum
$\nu_0 \in \Rgeq$.  Since $f \in \Lambda^0_{\{p\}}$,
we have $\nu_0 > 0$.  It is not hard to check that the
statement of the lemma holds for every $k > 1/\nu_0$.
\hfill$\Box$

\vspace{.125in}

Let $Z = X \cap (f^\perp)^{-1} (0)$.  Let $I^- = (-\infty, 0)$
and let $Z^- = Z \cap f^{-1} (I^-)$.  Consider the product
$\tilde Y^- = Y^- \times I^-$.  The decomposition of $Z^-$
into open rays $R^\circ(p,y)$, for $y \in Y^-$, defines an
isomorphism $\chi : \tilde Y^- \to Z^-$.  More precisely,
for $y \in Y^-$ and $a \in I^-$, we have
$$\{ \chi (y, a) \} = R^\circ (p,y) \cap f^{-1} (a).$$
Let $\tilde G^-$ be the vector field on $\tilde Y^-$ which
is given in components by $\tilde G^-_{(y,a)} = (G^-_y,0)$.
Also, consider the radial vector field
$E = - r \cdot  \nabla_\mu \, r |^{}_{Z^-}$ on $Z^-$.
Let $\pi : Z^- \to Y^-$ be the projection
$\pi : \chi (y, a) \mapsto y$.

\begin{lemma}\label{infprod}
There exist smooth functions $\alpha, \beta : Y^- \to \R_+$,
such that:
$$V|_{Z^-} = (\alpha \circ \pi) \cdot E +
(\beta \circ \pi) \cdot \chi_* (\tilde G^-).$$
\end{lemma}

\noindent
{\bf Proof:}  By $\R_+$-homogeneity, it is enough to find
$\alpha, \beta : Y^- \to \R_+$, such that
\begin{equation}\label{VEG}
V_y = \alpha (y) \cdot E_y + \beta (y) \cdot G^-_y,
\end{equation}
for every $y \in Y^-$.

Fix a point $y \in Y^-$.  Let $V^1_y = (\nabla_\mu f)^{}_y \in W$.
Let $T_y Y^- \subset T_y Z^- \subset W$ be the tangent spaces to
$Y^-$ and $Z^-$ at $y$, considered as a linear subspaces of $W$.
Note that
$$(\nabla_\mu \, r)^{}_y = - \frac{1}{r(y)} \cdot E_y \in T_y Z^-.$$
Let $\zeta : W \to T_y Z^-$ and $\eta : T_y Z^- \to T_y Y^-$ be
the orthogonal projections with respect to $\mu$.  By Lemma
\ref{imaginaryintegral}, we have $V_y = \zeta (V^1_y)$.  Also, we 
have  $G^-_y = \eta ((\nabla_\mu \, r)^{}_y)$.  Furthermore,
$T_y Y^- \subset T_y Z^-$ is the orthogonal complement to $V_y$.
Therefore, writing $\langle \; , \; \rangle$ for the real inner product
given by $\mu$, we have
$$G^-_y = -\frac{1}{r(y)} \cdot \eta (E_y) =
-\frac{1}{r(y)} \cdot \left( E_y - \frac{\langle E_y ,
V_y \rangle} {\langle V_y , V_y \rangle} \cdot V_y \right).$$
Note further that
$\langle E_y , V_y \rangle = \langle E_y , V^1_y \rangle =
E_y f = 1$.
Thus, equation (\ref{VEG}) holds with
$\alpha (y) = \langle V_y , V_y \rangle$ and
$\beta (y) =r(y) \cdot \langle V_y , V_y \rangle$.
It is easy to see that, so defined, the functions 
$\alpha, \beta : Y^- \to \R_+$ are smooth.
\hfill$\Box$

\vspace{.125in}

\noindent
{\bf Proof of Theorem \ref{model}:}
The claim of the theorem for the descending set $M^-_V (p)$
is proved by putting together Lemmas \ref{imaginaryintegral},
\ref{safeconical}, \ref{infprod}.

It is not hard to check that the flow of $V$ is integrable
for all time.  Let $\gamma : \R \to X^\circ$ be a parameterized
trajectory of $V$.  Let $\Gamma = \gamma (\R) \subset X^\circ$.
Suppose $\Gamma \subset M^-_V (p)$.  By Lemma
\ref{imaginaryintegral}, we have $\Gamma \subset Z^-$.
By Lemma \ref{infprod}, the image $\pi (\Gamma) \subset Y^-$
is a trajectory of $G^-$.  By Lemma \ref{safeconical}, the
trajectory $\pi (\Gamma)$ is bounded.  This proves the
containment $M^-_V (p) \subset \Cone (Y^-_c)$.  The opposite
containment is similar.

The claim of the theorem for the ascending set $M^+_V (p)$
is analogous.
\hfill$\Box$

\section{Proofs of the Main Results}

In this section, we will prove the following enhanced
versions of Theorems \ref{geom} and \ref{geom2}.

\begin{thm}\label{geomfadapted}
Theorem \ref{geom} is true and, moreover, the vector field
$V$ can be chosen to be compatible with a system of control
data which is $f$-adapted near $\{ p \}$.
\end{thm}

\begin{thm}\label{geom2fadapted}
Theorem \ref{geom2} is true and, moreover, the vector field
$V$ can be chosen to be compatible with a system of control
data which is $f$-adapted near $\{ p \}$.
\end{thm}

\subsection{Reduction to a Point Stratum}

The following proposition is a direct consequence of Corollaries
\ref{fadaptedexist}, \ref{fadaptedprod} and Proposition
\ref{fadaptednsequiv}.

\begin{prop}\label{reducetopoint}
Let $n > 0$ be an integer.  Assume that Theorem \ref{geomfadapted}
is true for $\dim_\C X < n$.  Then Theorem \ref{geom2fadapted} is
true for $\dim_\C X - \dim_\C A < n$.
\end{prop}

\noindent
{\bf Proof:}  Let $(X,\cS)$, $f: X \to \R$, $p \in \Sigma_f$, and
$A \in \cS$ be as in Theorem \ref{geom2}.  Assume that
$\dim_\C X - \dim_\C A < n$.  Use Corollary \ref{fadaptedexist} to
obtain a system of control data $\cD = \{ U_S, \Pi_S, \rho_S \}$ on
$(X, \cS)$ which is $f$-adapted near $\{ p \}$.  Let $N = \Pi_A^{-1}
(p)$; it is a normal slice to $A$.  Let $\cS_1$ be the stratification of
$N$ induced from $\cS$, let $\cD_1$ be the system of control data
on $(N, \cS_1)$ induced from $\cD$, and let $f_1 = f|^{}_N : N \to \R$.
Apply Corollary \ref{fadaptedprod} to obtain an open neighborhood
$\cU \subset U_A$ of $p$ and a controlled homeomorphism
$\phi : (\cU \cap N) \times (\cU \cap A) \to \cU$.

The normal slice $N$ need not be complex analytic, so we can not
apply Theorem \ref{geomfadapted} directly to the function $f_1 : N
\to \R$.  However, we can pick a complex analytic normal slice
$N_a$ to $A$ through $p$, apply Theorem \ref{geomfadapted} to
the function $f_a = f|^{}_{N_a} : N_a \to \R$, then use Proposition
\ref{fadaptednsequiv} to relate $f_1$ and $f_a$.  In this way, we
obtain an open neighborhood $\cU_1 \subset N$ of $p$ and a
$\nabla f_1$-like vector field $V_1$ on $\cU_1$, compatible with
$\cD_1$ and satisfying the conditions of Theorem \ref{geom}.
In particular, we obtain presentations of $L^\pm_{V_1} (p)$ and
$M^\pm_{V_1} (p)$ as cellular subsets of $\cU_1$.

By shrinking the neighborhoods $\cU$ and $\cU_1$, we can assume 
that $\cU_1 = \cU \cap N$.  Let $\cU_2 = \cU \cap A$, and let
$f_2 = f|^{}_{\cU_2} : \cU_2 \to \R$.  Fix a Riemannian metric
$\mu$ on $A$, and define a vector field $V_2$ on $\cU_2$ by
$V_2 = \nabla_\mu \, f_2$.  The ascending and descending links
$L^\pm_{V_2} (p)$ are both spheres, with
$$\dim_\R L^-_{V_2} (p) = \index_{f_2} (p) - 1, \;\; \mbox{and}$$
$$\dim_\R L^+_{V_2} (p) = \dim_\R A - \index_{f_2} (p) - 1.$$
(By convention, a sphere of dimension $-1$ is the empty set.)

Let us focus on the descending link $L^-_{V_2} (p)$.  Assuming
$L^-_{V_2} (p) \neq \emptyset$, fix a point $q^- \in L^-_{V_2} (p)$. 
Let $D^-_1 = \{ q^- \}$ and $D^-_2 = L^-_{V_2} (p) \setminus D^-_1$.
The following equation presents $L^-_{V_2} (p)$ as a cellular subset
of $\cU_2$:
$$L^-_{V_2} (p) = D^-_1 \cup D^-_2.$$
For $i = 1,2$, let $(D^-_i)^\sharp \subset M^-_{V_2} (p)$ be the union
of all trajectories of $V_2$ passing through $D^-_i$.  By shrinking
the neighborhood $\cU_2$, if necessary, we can assume that each
$(D^-_i)^\sharp$ is diffeomorphic to an open ball.  The following
equation then presents the descending set $M^-_{V_2} (p)$ as a
cellular subset of $\cU_2$:
$$M^-_{V_2} (p) = \{p\} \cup (D^-_1)^\sharp \cup (D^-_2)^\sharp.$$ 
In the case $L^-_{V_2} (p) = \emptyset$, we have $M^-_{V_2} (p) =
\{ p \}$, which is also a cellular subset of $\cU_2$.  We proceed
by analogy to present $L^+_{V_2} (p)$ and $M^+_{V_2} (p)$ as
cellular subsets of $\cU_2$.

We define $V = \phi_*(V^\times)$, where $V^\times$ is the vector
field on $\cU_1 \times \cU_2$ which is given in components by
$V_1$ and $V_2$.  Clearly, the vector field $V$ is compatible with
the system of control data $\cD$.  We have:
$$M^\pm_V (p) = \phi (M^\pm_{V_1} (p) \times M^\pm_{V_2} (p)).$$
A product of cellular subsets can be given a structure of a cellular
subset by ordering the pair-wise products of cells lexicographically.
This procedure presents the products
$$M^\pm_{V_1} (p) \times M^\pm_{V_2} (p) \subset 
\cU_1 \times \cU_2$$
as cellular subsets.  By applying the homeomorphism $\phi$,
we obtain structures of cellular subsets on the ascending and 
descending sets $M^\pm_V (p) \subset \cU$.  Verification of
conditions (ii)-(iv) of Theorem \ref{geom2} is routine.
\hfill$\Box$

\vspace{.125in}

Proposition \ref{reducetopoint} shows that Theorem
\ref{geomfadapted} implies Theorem \ref{geom2fadapted}.  We are
going to prove Theorem \ref{geomfadapted} by induction on
$d = \dim_\C X$.  The case $d = 0$ is trivial.  Fix an integer
$n > 0$, and assume that Theorem \ref{geomfadapted} is true for
$d < n$.  We will now proceed to establish Theorem
\ref{geomfadapted} for $d = n$, using Proposition
\ref{reducetopoint} in the process.

\subsection{Complex Links and Control Data}

Without loss of generality we can assume that $X$ is a
neighborhood of the origin in a complex vector space
$W \cong \C^d$ and that $p \in W$ is the origin.  Furthermore,
by Proposition \ref{fadaptednsequiv}, we can assume that
$f = \mbox{Re} (\varphi) |_X$, where $\varphi : W \to \C$
is a complex linear function.

Pick a Hermitian metric $\mu$ on $W$ and define $r : W \to \R$
by $r(w) = \mbox{dist}_\mu(p,w)$.  Let $\Delta = \mbox{Ker}
(\varphi) \subset W$, let $\theta : W \to \Delta$ be the orthogonal
projection with respect to $\mu$, and let $L = \theta^{-1} (p)
\subset W$.  Define $\hat r = r \circ \theta : W \to \R$.  Let
$$\Psi = (\varphi, \hat r): W \to \R^3 \; (= \C \times \R).$$
Define $\alpha : W \setminus L \to \R$ by
$\alpha (w) = |\varphi(w)|/\hat r(w)$.  Note that
$\hat r$ and  $\Psi$ are smooth on $W \setminus L$, while
$\alpha$ is smooth on $W \setminus (L \cup \Delta)$.

\begin{lemma}\label{saferegion}
There exist $\delta, \kappa > 0$ such that the following conditions
hold. 

(i)   Let $B_{\delta} = \{ w \in W \; | \; r(w) < \delta \}$.  Then
we have $B_{\delta} \subset X$.

(ii)  The map $\varphi$ is a stratified submersion on the
punctured ball $B_{\delta} \setminus \{ p \}$.

(iii) The map $\Psi$ is a stratified submersion on the region
$\{ x \in B_{\delta} \; | \; \alpha(x) < \kappa \}$.
\end{lemma}

\noindent
{\bf Proof:}  This is similar to the proof of \cite[Lemma 3.7]{Gr}.
\hfill$\Box$

\vspace{.125in}

We fix the numbers $\delta, \kappa > 0$ provided by Lemma
\ref{saferegion} for the rest of Section 4.  Let $f^\perp =
\mbox{Im} (\varphi) |_X : X \to \R$.  For $\epsilon > 0$, define
$$\cU_\epsilon = \{ x \in X \;\; | \;\;
|f(x)| < \epsilon,  \; |f^\perp (x)| < \epsilon, \;
\hat r(x) < 10 \cdot \epsilon / \kappa \}.$$
Fix an $\epsilon_0 > 0$ such that $\cU_{4\epsilon_0} \subset
B_\delta$.  We define the open neighborhood $\, \cU$ of
Theorems \ref{geom}, \ref{geomfadapted} by $\, \cU =
\cU_{2\epsilon_0}$.  Also, with a view to constructing a
system of control data on $(X, \cS)$, we let $U_{\{p\}} =
\cU_{3\epsilon_0}$.

Consider the manifolds
$$Y^\pm = \{ x \in \cU \;\; | \;\;
\varphi(x) = \pm \epsilon_0, \;
\hat r (x) < 6 \cdot \epsilon_0 / \kappa \},$$
and the restrictions
$$g^\pm = r|^{}_{Y^\pm}: Y^\pm \to \R.$$
We will be referring to the spaces $Y^\pm$ as the complex links
of $p$.  This is slightly different from standard terminology
(see \cite[Introduction, \S 1.5]{GM2}) in that $Y^\pm$ are not
closed in $X$.  By condition (ii) of Lemma \ref{saferegion},
the manifolds $Y^\pm$ meet the strata of $\cS$ transversely.
Therefore, each $Y^\pm$ is Whitney stratified by the
intersections with the strata of $\cS$.  We denote these
stratifications by $\cS^\pm$.  Consider the stratified
critical loci $\Sigma_{g^\pm} \subset Y^\pm$.  By condition
(iii) of Lemma \ref{saferegion}, we have
\begin{equation}\label{sigmagbound}
\Sigma_{g^\pm} \subset \{ y \in Y^\pm \; | \;
\hat r (y)  \leq \epsilon_0/\kappa \}.
\end{equation}
By perturbing the metric $\mu$ and reducing the numbers
$\delta, \kappa, \epsilon_0 >0$, if necessary, we may
assume that $g^\pm$ are stratified Morse functions
(this is similar to Proposition \ref{genmetric}).

Let $X^\circ = X \setminus \{ p \}$, let
$\cS^\circ = \{ S \in \cS \; | \; S \neq \{ p \} \}$
be the stratification of $X^\circ$ induced from $\cS$,
and let $U_{\{p\}}^\circ = U_{\{p\}} \setminus \{ p \}$.
Also, let $U_{\{p\}}^{\circ,1} = \{ x \in U_{\{p\}}^\circ \; |
\; \alpha(x) < \kappa/2 \}$.

\begin{lemma}\label{cdpunctured}
There exists a system of control data $\cD^\circ =
\{ U_S, \Pi_S, \rho_S\}_{S \in \cS^\circ}$ on $(X^\circ,
\cS^\circ)$, such that the following conditions hold.

(i)   The system $\cD^\circ$ is $\varphi$-compatible on
      $U_{\{p\}}^\circ$.  Moreover, for every $S \in \cS^\circ$,
      we have:
      $$\varphi \circ \Pi_S = \varphi \;\;\; \mbox{on} \;\;\;
      U_S \cap U_{\{p\}}^\circ$$
      (cf. Definition \ref{compat}). 
      
(ii)  The system $\cD^\circ$ is $\Psi$-compatible on
      $U_{\{p\}}^{\circ,1}$.  Moreover, for every $S \in \cS^\circ$,
      we have:
      $$\Psi \circ \Pi_S = \Psi \;\;\; \mbox{on} \;\;\;
      U_S \cap U_{\{p\}}^{\circ,1}.$$

(iii) Conditions (i)-(ii) imply that $\cD^\circ$ restricts to
      systems of control data $\cD^\pm$ on $(Y^\pm, \cS^\pm)$.
      For every $y \in \Sigma_{g^-}$, the system $\cD^-$ is
      $g^-$-adapted near $\{ y \}$; and similarly for $\cD^+$.
\end{lemma}

\noindent
{\bf Proof:}  This follows from Corollary \ref{fadaptedexist}
and Lemma \ref{saferegion} (cf. Step 4 of the proof of
\cite[Theorem 3.1]{Gr}).
\hfill$\Box$

\vspace{.125in}

Define $\rho_{\{p\}} : U_{\{p\}}^\circ \to \R_+$ by
\begin{equation}\label{definerho}
\rho_{\{p\}} (x) = \smax(
\smax((\kappa/10) \cdot \hat r(x),
\, |f^\perp (x)|), \, |f(x)|).
\end{equation}
Let $\Pi_{\{p\}} : U_{\{p\}} \to \{p\}$ be the unique
such map.

\begin{lemma}\label{cdfull}
Combining the system $\cD^\circ$ of Lemma \ref{cdpunctured}
with the triple $(U_{\{p\}}, \Pi_{\{p\}}, \rho_{\{p\}})$ produces a
system of control data $\cD = \{ U_S, \Pi_S, \rho_S\}_{S \in \cS}$
on $(X, \cS)$ which is $f$-adapted near $\{p\}$.
\end{lemma}

\noindent
{\bf Proof:}  Given that $\cD^\circ$ is a system of control data
on $(X^\circ, \cS^\circ)$, in order to check that $\cD$ is a system
of control data on $(X, \cS)$, we only need to verify conditions
(1)-(4) of Definition \ref{cdat} for $S = \{ p \}$.  Conditions (1)
and (2) follow from Lemma \ref{saferegion}.  Condition (3) is
essentially vacuous.  And condition (4) follows form conditions
(i)-(ii) of Lemma \ref{cdpunctured} and the fact that $\rho_{\{p\}}$
factors through $\Psi$ on $U_{\{p\}}^{\circ,1}$ and through $\phi$
on $U^\circ_{\{p\}} \setminus U_{\{p\}}^{\circ,1}$.  Verifying that
$\cD$ is $f$-adapted near $\{p\}$, based on Lemmas
\ref{saferegion} and \ref{cdpunctured}, is routine.
\hfill$\Box$

\subsection{Vector Fields on the Complex Links}

The following proposition will provide the key inductive input
into our construction of the vector field $V$ of Theorems
\ref{geom}, \ref{geomfadapted} in Section 4.4.  Define
$$Y^{\pm,s} = \{y \in Y^\pm \; | \;
\hat r (y) > 3 \cdot \epsilon_0 / \kappa \}.$$

\begin{prop}\label{clflow}
There exists a $\nabla g^-$-like vector field $G^-$ on
$Y^-$, compatible with the system of control data $\cD^-$
of Lemma \ref{cdpunctured}, such that the following
conditions hold.  

(i)  For every $y \in \Sigma_{g^-}$, there exists an open
neighborhood $\, \cU_y \subset Y^-$ of $y$ such that the
restriction $G^- |_{\cU_y}$ satisfies conditions (i)-(iv)
of Theorem \ref{geom2} with $(X, \cS) = (Y^-, \cS^-)$,
$f = g^-$, $p = y$, $\cU = \cU_y$, and $V = G^- |_{\cU_y}$.

(ii)  For every $y \in Y^{-,s}$, we have $G^-_y \, \hat r =
\hat r (y)$.

\noindent
Similarly, there exists a $\nabla g^+$-like vector field $G^+$
on $Y^+$ compatible with $\cD^+$ and satisfying the
analogues of conditions (i)-(ii) above.
\end{prop}

\noindent
{\bf Proof:}  We present the proof for $Y^-$ only.  Pick a critical
point $y \in \Sigma_{g^-}$.  By the induction hypothesis stated
at the end of Section 4.1 and by Proposition \ref{reducetopoint},
Theorem \ref{geom2fadapted} holds for the function $g^- :
Y^- \to \R$ and the critical point $y$.  This gives us an 
open neighborhood $\, \cU(y) \subset Y^-$ of $y$ and an
$\cS^-$-preserving $\nabla g^-$-like vector field $G(y)$ on
$\, \cU(y)$ satisfying conditions (i)-(iv) of Theorem \ref{geom2}.
Moreover, the vector field $G(y)$ can be chosen to be
compatible with a system of control data $\cD(y)$ on
$(Y^-, \cS^-)$ which is $g^-$-adapted near $\{ y \}$.
By Corollary \ref{fadaptedprod}, Proposition 
\ref{fadaptednsequiv}, and condition (iii) of Lemma
\ref{cdpunctured}, we can assume that $\cD(y) = \cD^-$.

Next, let $Y^{-,1} = \{ y \in Y^- \; | \; \hat r (y) > 2 \cdot
\epsilon_0 / \kappa \}$.  By equation (\ref{sigmagbound}),
we have $\Sigma_{g^-} \cap Y^{-,1} = \emptyset$.  By
condition (ii) of Lemma \ref{cdpunctured}, the system of
control data $\cD^-$ is $g^-$-compatible on $Y^{-,1}$.
Note that the function $\hat r$ factors through $g^-$ on
$Y^{-,1}$.  More precisely, let $I_1 = g^- (Y^{-,1})$ and
$I_2 = \hat r (Y^{-,1})$.  Then both $I_1, I_2 \subset \R$
are open intervals, and there exists an
orientation-preserving diffeomorphism $h: I_1 \to I_2$,
such that $\hat r = h \circ g^-$ on $Y^{-,1}$.
It follows that there exists a $\nabla g^-$-like vector field
$G^{-,1}$ on  $Y^{-,1}$ compatible with $\cD^-$, such that
$G^{-,1}_y \; \hat r = \hat r (y)$ for every $y \in Y^{-,s}$.

The requisite vector field $G^-$ is obtained by ``patching
together'' the vector fields $\{ G(y) \}$ for $y \in \Sigma_{g^-}$
and $G^{-,1}$, using Lemma \ref{ecvf} and a suitable
partition of unity on $Y^-$.  The only nuance is that $\cD^-$
must be compatible with the elements of the partition, to
ensure that $G^-$ satisfies condition (a) of Definition
\ref{wcvf}.
\hfill$\Box$

\subsection{Construction of the Vector Field $V$}

We are now prepared to construct the requisite $\nabla f$-like
vector field $V$ on the open set $\cU$ defined in Section 4.2. 
The vector field $V$ will be compatible with the system of
control data $\cD$ of Lemma \ref{cdfull}.  Let $Z = \cU \cap
(f^\perp)^{-1} (0)$.  By condition (ii) of Lemma \ref{saferegion},
the set $Z$ is a smooth manifold with a Whitney stratification
$\cZ$ induced from $\cS$.  By conditions (i)-(ii) of Lemma
\ref{cdpunctured}, the system of control data $\cD$ on
$(X, \cS)$ restricts to a system of control data $\check \cD$
on  $(Z, \cZ)$.  Let $\check f = f|_Z : Z \to \R$.  Define
$\cU^\circ = \cU \setminus \{ p \}$, $Z^\circ = Z \setminus
\{ p \}$, and write $\rho = \rho_{\{p\}} : U_{\{p\}}^\circ \to \R_+$
to unclutter the notation.

\begin{lemma}\label{extendfromz}
Let $\check V$ be a $\nabla \check f$-like vector field on
$Z$ compatible with $\check \cD$.  Then there exists a
$\nabla f$-like vector field $\hat V$ on $\, \cU$ compatible
with $\cD$, such that:

(i)  $\hat V |_Z = \check V$;

(ii) $\hat V_x \, f^\perp = 0$ for every $x \in \cU$.

\noindent
Moreover, for every such $\hat V$, we have
$M^\pm_{\hat V} (p) = M^\pm_{\check V} (p) \subset Z$.
\end{lemma}

\noindent
{\bf Proof:}
Cover the set $\cU^\circ$ by two open subsets as follows:  
$$\cU^{\circ,1} = \{
x \in \cU^\circ \; | \; | f^\perp (x) | < \rho(x) / 10 \},$$
$$\cU^{\circ,2} = \{
x \in \cU^\circ \; | \; | f^\perp (x) | > \rho(x) / 20 \}.$$
Let $I = (-2 \epsilon_0, 2 \epsilon_0) \subset \R$.
Define an open subset $\tilde Z^\circ \subset Z^\circ
\times I$ by
$$\tilde Z^\circ = \{ (x, a) \in Z^\circ \times I \; | \;
|a| < \rho(x) / 10 \}.$$
Consider $\tilde Z^\circ$ as a smooth manifold with a
Whitney stratification $\tilde \cZ^\circ$ induced from
$\cZ$ and a system of control data $\tilde \cD^\circ$
induced from $\check \cD$.  By Lemmas \ref{ecvf},
\ref{saferegion}, \ref{cdpunctured} and equation
(\ref{definerho}), there exists a controlled
homeomorphism
$$h : \tilde Z^\circ \to \cU^{\circ,1},$$
compatible with $\tilde \cD^\circ$ and $\cD$, such that
for every $(x, a) \in \tilde Z^\circ$, we have:

(a)  $h \; (x,0) = x$;

(b)  $f^\perp \circ h \; (x,a) = a$;

(c)  $f \circ h \; (x,a) = f(x)$;

(d)  $\rho \circ h \; (x,a) = \rho(x)$.

\noindent
Define a vector field $\hat V^1$ on $\cU^{\circ,1}$ by
$\hat V^1_{(x,a)} = h_*(\check V_x,0)$.  Use Lemmas 
\ref{ecvf}, \ref{saferegion}, \ref{cdpunctured} and
equation (\ref{definerho}) again, to construct a
$\nabla f$-like vector field $\hat V^2$ on $\cU^{\circ,2}$
compatible with $\cD$, such that for some $k_2 > 0$ and every
$x \in \cU^{\circ,2}$, we have $\hat V^2_x \, f^\perp = 0$ and
$| \hat V^2_x \, \rho | < k_2 \cdot \rho(x)$.  The restriction
of the requisite vector field $\hat V$ to $\cU^\circ$ is
constructed by combining the vector fields $\hat V^1$ and
$\hat V^2$, using a suitable partition of unity on $\cU^\circ$.
Verification of conditions (i) and (ii) is routine.  The claim
of the lemma about $M^\pm_{\hat V} (p)$ follows immediately
from conditions (i) and (ii).
\hfill$\Box$

\vspace{.125in}

We will first construct the restriction $\check V = V|_Z$,
then use Lemma \ref{extendfromz} to obtain the full vector
field $V$.  Define
$$Z^- = \{ x \in Z^\circ \;\; | \;\;
f(x) < 0, \; \hat r (x) < (6/\kappa) \cdot |f(x)| \},$$
$$Z^s = \{ x \in Z^\circ \;\; | \;\; 
\hat r (x) > (3/\kappa) \cdot |f(x)| \},$$
$$Z^{ss} = \{ x \in Z^\circ \;\; | \;\; 
\hat r (x) > (4/\kappa) \cdot |f(x)| \},$$
$$Z^+ = \{ x \in Z^\circ \;\; | \;\;
f(x) > 0, \; \hat r (x) < (6/\kappa) \cdot |f(x)| \}.$$
The superscript ``$s$'' stand for ``safe.''  In the construction
that follows, any trajectory of $\check V$ which enters $Z^s$
will be safe from approaching the critical point $p$. 
We begin by constructing vector fields $\check V^\pm$ on
$Z^\pm$.  They will serve as the main building blocks in the
construction of $\check V$ (see Lemma \ref{vwaist}).  We
will only describe the construction of $\check V^-$.  The
construction of $\check V^+$ is analogous.  Lemmas
\ref{manmadeprodvf}, \ref{manmadeprod}, \ref{vwaist},
\ref{extendfromznot} below should be seen as
Lemmas/Definitions; they introduce objects which will
be referred to directly later.

\begin{lemma}\label{manmadeprodvf} 
There exists a controlled vector field $E$ on $Z^-$
compatible with $\check \cD$, such that the following
conditions hold.

(i)  For every $x \in Z^-$, the derivative $E_x \, f =
     - f(x)$.

(ii) For every $x \in Z^- \cap Z^s$, the derivative
$E_x \, \hat r = - \hat r (x)$.
\end{lemma}

\noindent
{\bf Proof:}  This is similar to the proof of Proposition
\ref{clflow}, using Lemmas \ref{ecvf}, \ref{saferegion},
\ref{cdpunctured} and a partition of unity argument.
\hfill$\Box$

\vspace{.125in}

Let $I^- = (-2\epsilon_0,0)$. Consider the product $\tilde Y^- =
Y^- \times I^-$.  Let $\pi_1 : \tilde Y^- \to Y^-$ and $\pi_2 :
\tilde Y^- \to I^-$ be the projection maps.  Let $\tilde \cS^-$
be the stratification of $\tilde Y^-$ induced from $\cS^-$, and
let $\tilde \cD^-$ be the system of control data on $(\tilde Y^-,
\tilde \cS^-)$ induced from $\cD^-$.

\begin{lemma}\label{manmadeprod} 
There exists a unique controlled homeomorphism $\chi :
\tilde Y^- \to Z^-$ compatible with $\tilde \cD^-$ and
$\check \cD$, such that the following conditions hold.

(i)   For every $y \in Y^-$ and every $\epsilon \in I^-$,
      we have $f(\chi(y,\epsilon)) = \epsilon$.

(ii)  For every $y \in Y^-$, we have $\chi(y, -\epsilon_0) = y$.

(iii) For every $y \in Y^-$, the set $\chi(\pi_1^{-1} (y))$
      is a trajectory of $E$.

(iv)  For every $y \in Y^- \cap Z^s$ and every
      $\epsilon \in I^-$, we have
      $\alpha(\chi(y,\epsilon)) = \alpha(y)$.
\end{lemma}

\noindent
{\bf Proof:}  It is not hard to check that the requisite
homeomorphism $\chi$ is given by:
$$\chi (y, \epsilon) = \psi_{E,t(\epsilon)} (y),$$
where $t(\epsilon) = \ln(\epsilon_0/\epsilon)$ and 
$\psi_{E,t(\epsilon)}$ is the flow of the vector field $E$.
\hfill$\Box$

\vspace{.125in}

Recall the vector field $G^-$ on $Y^-$, provided by
Proposition \ref{clflow}.  Let $\tilde G^-$ be the weakly
controlled vector field on $\tilde Y^-$ which is given in
components by $\tilde G^-_{(y,\epsilon)} = (G^-_y,0)$.
We define
\begin{equation}\label{infprod2}
\check V^-= E + \chi_*(\tilde G^-).
\end{equation}

\begin{lemma}\label{vminus}
The vector field $\check V^-$ is $\nabla \check f$-like
on $Z^-$.  Furthermore, it satisfies the following conditions.

(i)  For every $x \in Z^-$, we have $\check V^-_x \, f =
     \rho (x)$.

(ii) For every $x \in Z^- \cap Z^s$, we have
     $\check V^-_x \, \hat r = 0$.
\end{lemma}

\noindent
{\bf Proof:}  By condition (i) of Lemma \ref{manmadeprod},
we have $\chi_*(\tilde G^-)_x \, f = 0$ for every $x \in Z^-$.
Combining this with condition (i) of Lemma
\ref{manmadeprodvf}, we obtain:
$$\check V^-_x \, f = E_x \, f = - f(x) = |f(x)|.$$
By equation (\ref{definerho}), we have $\rho(x) = |f(x)|$ for
every $x \in Z^-$.  This verifies condition (i) and the fact
that $\check V^-$ is $\nabla \check f$-like.

Next, by condition (ii) of Proposition \ref{clflow} and
conditions (i), (iv) of Lemma \ref{manmadeprod}, we
have $\chi_*(\tilde G^-)_x \, \hat r = \hat r (x)$ for every
$x \in Z^- \cap Z^s$.  Combining this with condition (ii) of
Lemma \ref{manmadeprodvf}, we obtain:
$$\check V^-_x \, \hat r = \chi_*(\tilde G^-)_x \, \hat r
+ E_x \, \hat r = \hat r (x) - \hat r (x) = 0.$$
This verifies condition (ii).
\hfill$\Box$

\vspace{.125in}

At this point, we assume that we have carried out the analogous
construction to obtain a $\nabla \check f$-like vector field
$\check V^+$ on $Z^+$, satisfying $\check V^+_x f = \rho (x)$
for every $x \in Z^+$ and $\check V^+_x \, \hat r = 0$ for every
$x \in Z^+ \cap Z^s$.

\begin{lemma}\label{vwaist}
There exists a $\nabla \check f$-like vector field $\check V^\circ$
on $Z^\circ$ compatible with $\check \cD$, such that the following
conditions hold.

(i)   For every $x \in Z^\pm \setminus Z^{ss}$, we have
      $\check V^\circ_x = \check V^\pm_x$.

(ii)  For every $x \in Z^\circ$, we have
      $\check V^\circ_x \, f = \rho (x)$.

(iii) For every $x \in Z^s$, we have
      $\check V^\circ_x \, \hat r = 0$.
\end{lemma}

\noindent
{\bf Proof:}  This is another partition of unity argument,
similar to the proof of Proposition \ref{clflow}.  By Lemma
\ref{ecvf}, condition (iii) of Lemma \ref{saferegion}, condition
(ii) of Lemma \ref{cdpunctured}, and equation (\ref{definerho}),
there exists a weakly controlled vector field
$V^{ss}$ on $Z^{ss}$ compatible $\check \cD$, such that
$\check V^{ss}_x \, f = \rho (x)$ and 
$\check V^{ss}_x \, \hat r = 0$ for every $x \in Z^{ss}$.  Pick
a partition of unity $\eta^- + \eta^{ss} + \eta^+ = 1$ on
$Z^\circ$, such that $\supp (\eta^\pm) \subset Z^\pm$,
$\supp (\eta^{ss}) \subset Z^{ss}$, and each of the functions
$\eta^\pm, \eta^{ss} : Z^\circ \to [0,1]$ factors through
$\alpha : W \setminus L \to \R$ on the overlaps
$Z^\pm \cap Z^{ss} \subset W \setminus L$.  We define
$\check V^\circ$ by
$$\check V^\circ_x = \eta^- (x) \cdot \check V^-_x
+ \eta^{ss} (x) \cdot \check V^{ss}_x
+ \eta^+ (x) \cdot \check V^+_x,$$
for every $x \in Z^\circ$.  In the above equation, each of the
terms $\eta^* (x) \cdot \check V^*_x$ is understood to be zero
if $\eta^* (x) = 0$ and $\check V^*_x$ is undefined.
Verification of the properties of $\check V^\circ$ is routine.
\hfill$\Box$

\begin{lemma}\label{extendfromznot}
Define a vector field $\check V$ on $Z$ by $\check V|_{Z^\circ} =
\check V^\circ$ and $\check V_p = 0$.  Then $\check V$ is a
$\nabla \check f$-like vector field on $Z$ compatible with
$\check \cD$.
\end{lemma}

\noindent
{\bf Proof:}  By Lemma \ref{vwaist}, the restriction
$\check V|_{Z^\circ}$ is a $\nabla \check f$-like vector
field compatible with $\check \cD$.  The only property of
$\check V$ left to verify is condition (b) of Definition
\ref{wcvf} for the stratum $\{ p \} \in \cZ$.  We will show
that there exists a $k > 0$, such that
\begin{equation}\label{condb}
|\check V_x \, \rho| < k \cdot \rho (x),
\end{equation}
for every $x \in Z^\circ$.  Recall the function $\smax :
\Rgeq \times \Rgeq \to \Rgeq$ introduced in Section 3.3
and used in equation (\ref{definerho}).  By property (1) of
$\smax$, we can consider the partial derivative
$$\zeta(x,y) = \frac{\partial \, \smax(x,y)}{\partial \, x},$$
as a function $\zeta : \R_+ \times \R_+ \to \R$.  By property
(4) of $\smax$, we have $\zeta (a \cdot x, a \cdot y) =
\zeta (x,y)$ for all $a,x,y \in \R_+$.  By property (2) of
$\smax$, we have $\zeta(x,y) = 1$ for all  $x,y \in \R_+$
with $x > 1.1 \cdot y$, and $\zeta(x,y) = 0$ for all $x,y
\in \R_+$ with $y > 1.1 \cdot x$.  It follows that the function
$\zeta$ attains a maximum $c = \max(\zeta) \in \R$.
Conditions (ii) and (iii) of Lemma \ref{vwaist} then imply
that inequality (\ref{condb}) holds for every $k > c$.
It is a pleasant exercise, not necessary for this proof,
to check that, in fact, $c = 1$.
\hfill$\Box$

\vspace{.125in}

We use Lemma \ref{extendfromz} to obtain an extension $\hat V$
of $\check V$ to $\cU$; then let $V = \hat V$.  This completes
our construction of the vector field $V$.

\subsection{Identification of the Sets $M^\pm_V (p)$}

In this section, we identify the sets $M^\pm_V (p)$ for the
vector field $V$ constructed in Section 4.4, and verify
conditions (i)-(iv) of Theorem \ref{geom}.  We will limit our
discussion to the descending set $M^-_V (p)$.  The discussion
for the ascending set $M^+_V (p)$ is analogous and will be
omitted.

Recall the vector field $G^-$ on the complex link $Y^-$,
provided by Proposition \ref{clflow}.  For every $y \in
\Sigma_{g^-}$, consider the descending set $M^-_{G^-} (y)
\subset Y^-$.  Define
\begin{equation}\label{ycdef}
Y^-_c = \bigcup_{y \in \Sigma_{g^-}} M^-_{G^-} (y).
\end{equation}
By the continuity of the flow of $G^-$, the set
$Y^-_c \subset Y^-$ is compact.  Moreover, by equation
(\ref{sigmagbound}), we have
$$Y^-_c \subset \{ y \in Y^- \; | \; 
\hat r (y) \leq \epsilon_0 / \kappa \}.$$
Recall the controlled homeomorphism $\chi : \tilde Y^- \to Z^-$
of Lemma \ref{manmadeprod}.  The following proposition is
analogous to Theorem \ref{model}.

\begin{prop}\label{identifyds}
We have $M^-_V (p) = \{ p \} \cup \chi(Y^-_c \times I^-)$.
\end{prop}
 
\noindent
{\bf Proof:}  This is analogous to the proof of Theorem
\ref{model}.  The role of Lemma \ref{imaginaryintegral}
in that proof is played by condition (ii) of Lemma
\ref{extendfromz}.  The role of Lemma \ref{safeconical}
is played by condition (iii) of Lemma \ref{vwaist}, which
implies that $M^-_V (p) \cap Z^s = \emptyset$.  The role
of Lemma \ref{infprod} is played by equation (\ref{infprod2}).
\hfill$\Box$

\vspace{.125in}

Proposition \ref{identifyds} asserts that $Y^-_c$ is a descending
link of $p$.  We will now show that $Y^-_c$ is naturally a 
cellular subset of $Y^-$.  As in the proof of Corollary
\ref{modelcells}, order the set $\Sigma_{g^-}$ by critical
value.  More precisely, fix an ordering $\Sigma_{g^-} =
\{ y_0, \dots, y_N \}$ satisfying implication (\ref{ordercrit})
from that proof.
For each $y \in \Sigma_{g^-}$, recall the open neighborhood
$\cU_y \subset Y^-$ of Proposition \ref{clflow}, pick a small
number $\nu < 0$, and consider the descending link
$L^-(y) = L^-_{G^-,\nu} (y) \subset \cU_y$.  By Proposition
\ref{clflow}, the set $L^-(y)$ is a cellular subset of $\cU_y$.
Let $\cL^-(y) = \{ C_0 (y), \dots, C_{n(y)} (y)\}$ be the 
ordered set of cells of $L^-(y)$.  For $C \in \cL^-(y)$, define
$C^\sharp \subset M^-_{G^-} (y)$ to be the union of all
trajectories $\Gamma$ of $G^-$ with $\Gamma \cap C \neq
\emptyset$.  It is not hard to check that $C^\sharp$ is
diffeomorphic to an open ball with $\dim_\R C^\sharp =
\dim_\R C + 1$.  Furthermore, we have
\begin{equation}\label{mycells}
M^-_{G^-} (y) = \{ y \} \; \cup \bigcup_{C \in \cL^-(y)} C^\sharp.
\end{equation}

\begin{prop}\label{celldecomp}
(i)  We have:
\begin{equation}\label{ordercells}
Y^-_c = \bigcup_{i = 0}^{N} \left( \{y_i\} \cup
\bigcup_{j = 0}^{n(y_i)} C_j(y_i)^\sharp  \right).
\end{equation}
Moreover, the above equation endows $Y^-_c \subset Y^-$
with a structure of a cellular subset.

(ii) For every stratum $S \in \cS^-$, we have
$\dim_\R Y^-_c \cap S \leq \dim_\C S$.
\end{prop}

\noindent
{\bf Proof:}  
Equation (\ref{ordercells}) follows from equations
(\ref{ycdef}) and (\ref{mycells}).  The right-hand side of
(\ref{ordercells}) is an ordered union of cells.  Condition
(b) of Definition \ref{wsss} follows from the fact each
$L^-(y_i) \subset Y^-$ is a cellular subset and the
continuity of the flow of $G^-$.  This verifies claim (i).

Claim (ii) follows from condition (i) of Proposition \ref{clflow} 
and the inequality $\index_{g^-} (y_i) \leq 0$ for every
$y_i \in \Sigma_{g^-}$ (cf. inequality (\ref{indexbound}) in
the proof of Corollary \ref{modelcells}).
\hfill$\Box$

\vspace{.125in}

We are now prepared to describe a cell decomposition of
$M^-_V (p)$.  For every $y \in \Sigma_{g^-}$ and $C \in
\cL^-(y)$, define $\{ y \}^\sharp = \chi(\{ y \} \times I^-)$
and $C^{\sharp\sharp} = \chi(C^\sharp \times I^-)$.

\begin{prop}\label{celldecomp2}
We have
$$M^-_V (p) = \{ p \} \; \cup \;
\bigcup_{i = 0}^{N} \left( \{y_i\}^\sharp \cup
\bigcup_{j = 0}^{n(y_i)} C_j(y_i)^{\sharp\sharp}  \right).$$
Moreover, the above equation endows $M^-_V (p)
\subset \cU$ with a structure of a cellular subset.
\end{prop}

\noindent
{\bf Proof:}  This follows from Proposition \ref{identifyds},
claim (i) of Proposition \ref{celldecomp}, and the continuity
of the flow of $V$.
\hfill$\Box$

\vspace{.125in}

Condition (i) of Theorem \ref{geom} for the vector field $V$
follows from Proposition \ref{celldecomp2}.  Condition (ii)
follows from equation (\ref{infprod2}) and condition (iii) of
Lemma \ref{manmadeprod}.  Conditions (iii) and (iv) follow
from claims (i)  and (ii) of Proposition \ref{celldecomp},
respectively.  Since the vector field $V$ is compatible with
the system of control data $\cD$ of Lemma \ref{cdfull},
which is $f$-adapted near $\{ p \}$, this completes our proof
of Theorem \ref{geomfadapted}.  Theorem \ref{geom2fadapted}
follows by Proposition \ref{reducetopoint}.

\vspace{.1in}

\noindent
34 Hamlet Woods Drive, St James, NY 11780\hfill
{\em misha@mishagrinberg.net}

\end{document}